\documentclass[10pt]{amsart}
\usepackage{latexsym}
\usepackage{amsmath}
\usepackage{amssymb}
\usepackage{amscd}
\usepackage{array}
\usepackage{hhline}
\usepackage{color}
\usepackage[all]{xy}
\usepackage{enumerate}
\usepackage{graphicx}
\usepackage{lscape}

\theoremstyle{break}
\newtheorem{de}{Definition}[section]
\newtheorem{thm}[de]{Theorem}
\newtheorem{pro}[de]{Proposition}
\newtheorem{rem}[de]{Remark}

\newtheorem{cor}[de]{Corollary}
\newtheorem{conj}[de]{Conjecture}

\newtheorem{str}[de]{Strategy}

\def\C{{\mathbb{C}}}
\def\Q{{\mathbb{Q}}}
\def\N{{\mathbb{N}}}
\def\Z{{\mathbb{Z}}}
\def\c{{\mathcal{C}}}
\def\deg{{\mathrm{deg}}}
\def\dim{{\mathrm{dim}}}
\def\mod{{\mathrm{mod\;}}}

\def\ob{{\mathrm{Ob}}}
\def\mor{{\mathrm{Mor}}}
\def\M{{\mathrm{Mat}}}
\def\H{{\mathrm{H}}}
\def\kom{{\mathrm{Kom}}}
\def\id{{\mathrm{Id}}}
\def\MF{{\mathrm{MF}}}
\def\HMF{{\mathrm{HMF}}}

\def\ostimes{{\,\otimes\hspace{-0.7em}\raisebox{-0.5ex}{${}_{{}_{S}}$}\,}}
\def\osatimes{{\,\otimes\hspace{-0.7em}\raisebox{-0.5ex}{${}_{{}_{S\acute{}}}$}\,}}

\def\oqtimes{{\,\otimes\hspace{-0.6em}\raisebox{-0.5ex}{${}_{{}_{\Q}}$}\,}}

\title[Matrix factorizations and $U_q (\mathfrak{sl}_n)$ intertwiners]{Matrix factorizations and intertwiners of the fundamental representations of quantum group $U_q (\mathfrak{sl}_n)$}
\author{Yasuyoshi Yonezawa}
\address{Graduate School of Mathematics, Nagoya University\\ 464-8602 Furocho, Chikusaku, Nagoya, Japan }
\email{yasuyoshi.yonezawa@math.nagoya-u.ac.jp}
\date{}
\begin{document}
\maketitle
\begin{abstract}
We want to construct a homological link invariant whose Euler characteristic is MOY polynomial as Khovanov and Rozansky constructed a categorification of HOMFLY polynomial. The present paper gives the first step to construct a categorification of MOY polynomial. For the essential colored planar diagrams with additional data which is a sequence naturally induced by coloring, we define matrix factorizations, and then we define a matrix factorization for planar diagram obtained by gluing the essential planar diagrams as tensor product of the matrix factorizations for the essential planar diagrams. Moreover, we show that some matrix factorizations derived from tensor product of the essential matrix factorizations have homotopy equivalences corresponding to MOY relations.
\end{abstract}
\tableofcontents

%%%%%%%%%%%%%%%%%%%%%%%%%%%%%%%%%%%%%%%%%%%%%%%%%%%%%%%%%%%%%%%%%%%%%%%%%%%%%%%%%%%%%%%%%%%%%%%%%%%%%%%%
%
%
% section 1 : Introduction
%
%
%%%%%%%%%%%%%%%%%%%%%%%%%%%%%%%%%%%%%%%%%%%%%%%%%%%%%%%%%%%%%%%%%%%%%%%%%%%%%%%%%%%%%%%%%%%%%%%%%%%%%%%%
\section{Introduction}
\subsection{Categorification of quantum link invariant}\label{1-cat}
Mikhail Khovanov introduced a homological link invariant whose Euler characteristic is Jones polynomial. At present we understand Jones polynomial as a link invariant derived from the quantum group $U_q(\mathfrak{sl}_2)$ and its $2$-dimensional vector representation $V_2$. (We can also obtain a link invariant induced by a quantum group and its representation, called a quantum link invariant.)
Such a system making a homological link invariant whose Euler characteristic is a quantum link invariant is called a categorification of the quantum link invariant.
A natural question is ``can we construct a categorification of the other quantum link invariants?''
In the case of HOMFLY polynomial, which is derived from $U_q(\frak{sl}_n)$ and its $n$-dimensional representation $V_n$, Mikhail Khovanov and Lev Rozansky also constructed a homological link invariant whose Euler characteristic is HOMFLY polynomial.
\begin{figure}[hbt]
\begin{eqnarray*}
\left<\input{figplus}\right>_n&=&q^{-1+n}\left<\input{figsmoothing1sln}\right>_n-q^{n}\left<\input{figsmoothing2sln1}\right>_n\\[-0.1em]
\left<\input{figminus}\right>_n&=&q^{1-n}\left<\input{figsmoothing1sln}\right>_n-q^{-n}\left<\input{figsmoothing2sln1}\right>_n
\end{eqnarray*}
\caption{Reductions for single crossings of HOMFLY polynomial}
\end{figure}
\\ 
%%%%%%%%%%%%%%%%%%%%%%%%%%%%%%%%%%%%%%%
\indent However, there exist a lot of quantum link invariants which is not categorified yet. For example, MOY polynomial (see \cite{MOY}), which is a quantum (regular) link invariant derived from $U_q(\frak{sl}_n)$ and its fundamental representations, is one of quantum link invariants which are not categorified (uncategorified) yet.
Since this MOY polynomial is a generalization of HOMFLY polynomial, it is natural that we hope to construct a categorification of MOY polynomial generalizing Khovanov and Rozansky's work. \\
\indent We briefly recall the work of M.Khovanov and L.Rozansky \cite{KR1}. We calculate HOMFLY polynomial of an oriented link diagram $D$ by transforming each single crossing into planar diagrams $P_0$ and $P_1$ in Figure $2$ (in the case of Hopf link, see Figure $3$), evaluating the planar diagrams $\Gamma$ derived from the link diagram as a Laurent polynomial of $q$ and summing the Laurent polynomials by the reduction in Figure $1$. \\
\begin{figure}[hbt]
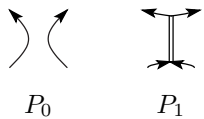

\begin{equation*}
\input{figsmoothing1sln-p}\hspace{1cm}\input{figsmoothing2sln1-p}
\end{equation*}
\caption{Planar diagrams}
\end{figure}
\begin{figure}[hbt]
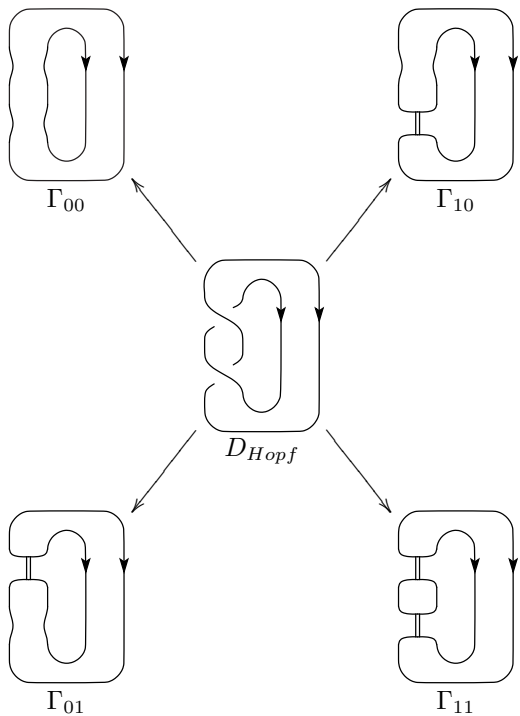

\begin{equation*}
\xymatrix{
\input{hopf-00}&&\input{hopf-10}\\
&\input{hopf-link}\ar[ru]\ar[rd]\ar[lu]\ar[ld]&\\
\input{hopf-01}&&\input{hopf-11}
}
\end{equation*}
\caption{Planar diagrams derived from Hopf link diagram}
\end{figure}
\indent For an oriented link diagram $D$, we can construct a homological link invariant $\c(D)_n$ whose Euler characteristic is HOMFLY polynomial for the link diagram $D$ as follows. 
First, for each planar diagram $P$ in Figure $2$ we define a matrix factorization $\c\acute{}(P)_n$ which is a 2-cyclic complex of modules.
Since every planar diagram $\Gamma$ induced by a link diagram $D$ is decomposed into some of the planar diagrams $P_0$ and $P_1$, a matrix factorization $\c\acute{}(\Gamma)_n$ for the planar diagrams $\Gamma$ is defined to be the tensor product of matrix factorizations for parts $P_0$ and $P_1$ of the decomposition. 
For example, the planar diagram $\Gamma_{10}$ in Figure $3$ has a decomposition shown as Figure $4$. Then we have a matrix factorization $\c\acute{}(\Gamma_{10})_n$ for the diagram $\Gamma_{10}$ as
\begin{equation*}
\c\acute{}(\Gamma_{10})_n= \c\acute{}(\Gamma_1)_n\boxtimes\c\acute{}(\Gamma_2)_n\boxtimes\c\acute{}(\Gamma_3)_n.
\end{equation*}
\begin{figure}[hbt]
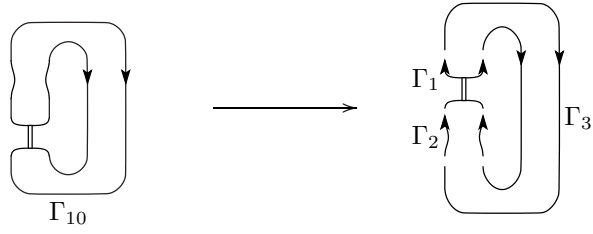

\begin{equation*}
\xymatrix{
\input{hopf-10}&\ar[rr]&&&\input{hopf-10-part}
}
\end{equation*}
\caption{Decomposition of planar diagram}
\end{figure}
\\
\indent For a link diagram $D$, we define a complex $\c(D)_n$ by exchanging every crossing into a complex of matrix factorizations for planar diagrams $P_0$ and $P_1$ as follows,
\begin{eqnarray*}
&&\xymatrix{
\c\Big( \input{figplus} \Big)_n  := \Big( \ar[r] 
&*^{\stackrel{-1\atop}{\c\acute{}\Big( \input{figsmoothing2sln1-mf} \Big)_n\{ n\} \left< 1 \right>}}\ar[r]^{\chi_+}
&*^{\stackrel{0\atop}{\c\acute{}\Big( \input{figsmoothing-sln-mf} \Big)_n\{ n-1 \} \left< 1 \right>}} \ar[r] 
&*^{0}\ar[r] & \Big), } \\[-0.5em]
&&
\xymatrix{\c\Big( \input{figminus} \Big)_n  := \Big(\ar[r] 
&*^{0}\ar[r] 
&*^{\stackrel{0\atop}{ \c\acute{}\Big( \input{figsmoothing-sln-mf} \Big)_n\{ -n+1 \} \left< 1 \right>}}\ar[r]^{\chi_-} 
&*^{\stackrel{1\atop}{ \c\acute{}\Big( \input{figsmoothing2sln1-mf} \Big)_n\{ -n \} \left< 1 \right>}}\ar[r] &\Big).
}
\end{eqnarray*}
\indent In the case of Hopf link diagram, we obtain the following complex,
\begin{equation*}
\xymatrix{
\c\acute{}(\Gamma_{00})_n\ar[rrr]^(0.4){\left(\begin{array}{c}{}_{\chi_-\boxtimes \id} \\ {}_{\id\boxtimes\chi_-} \end{array}\right)}&&&\c\acute{}(\Gamma_{10})_n\oplus\c\acute{}(\Gamma_{01})_n\ar[rrr]^(0.6){(\id\boxtimes\chi_-,-\chi_-\boxtimes\id)}&&&\c\acute{}(\Gamma_{11})_n
} .
\end{equation*}
\indent M. Khovanov and L. Rozansky introduced such a complex of matrix factorizations for a link diagram $D$ and they proved that $\chi_+$ and $\chi_-$ of the complex are associated with a homological link invariant, i.e. they showed that complexes derived from diagrams appearing in Reidemeister moves are isomorphic. \\
%%%%%%%%%%%%%%%%%%%%%%%%%%%%%%%%%%%%%%%
\indent We want to construct a homological link invariant whose Euler characteristic is MOY polynomial in \cite{MOY} as Khovanov and Rozansky constructed a categorification of HOMFLY polynomial. The present paper gives the first step to construct a categorification of MOY polynomial.
We consider the following strategy to achieve this purpose, 
\begin{str}[Categorification of MOY polynomial]\label{stra}
\indent\\
\begin{itemize}
\item[(S1)]Give a graded category $\MF^{gr}$ of a $\Z$-graded matrix factorization with tensor product (commutativity and associativity) such that some matrix factorization $\c\acute{}(P)$ of $\MF^{gr}$ realizes every essential colored planar diagram $P$ of Figure $5$.
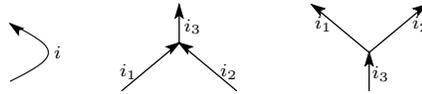
\begin{figure}[htb]
\input{figmoy3} \hspace{.5cm} \input{figgluing-in-3valent} \hspace{.5cm} \input{figgluing-out-3valent}\\
\begin{equation*}
1\leq i\leq n, 1\leq i_1,i_2\leq n-1, i_3=i_1+i_2\leq n
\end{equation*}
\caption{Essential colored planar diagrams}
\end{figure}
\\
\indent By calculation of MOY polynomial (see \cite{MOY}) every colored crossing is exchanged into a formal linear sum of $\Gamma^L_k$ and $\Gamma^R_k$ in Figure $6$. Since the diagrams $\Gamma^L_k$ and $\Gamma^R_k$ is decomposed into some of essential planar diagrams in Figure $5$, every colored planar diagram $\Gamma$ derived from a colored link diagram $D$ is also decomposed into some of essential planar diagrams in Figure $5$. For the colored planar diagram $\Gamma$, we define a matrix factorization $\c\acute{}(\Gamma)$ to be tensor product of matrix factorizations for parts of the decomposition.
\begin{figure}[hbt]
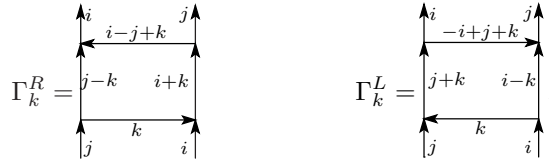
\begin{eqnarray*}
\Gamma^R_k=\input{figmoysmooth4}\hspace{1cm} \hspace{1cm} \Gamma^L_k=\input{figmoysmooth3},\\[1em]
\end{eqnarray*}
\caption{Colored planar diagrams}
\end{figure}
\\
Moreover, the category $\MF^{gr}$ must have a homotopy equivalence (i.e. an isomorphism in the homotopy category $\HMF^{gr}$) corresponding to every relation of MOY polynomial. For example, the homotopy category $\HMF^{gr}$ must have the following equivalence corresponding to MOY relation in \cite{MOY},
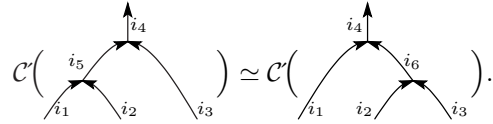
\begin{eqnarray*}
\c\acute{}\Big( \input{fig-ass-dia1}\Big) \simeq \c\acute{}\Big( \input{fig-ass-dia2}\Big) .\\[1em]
\end{eqnarray*}
\item[(S2)]In such a homotopy category $\HMF^{gr}$, we construct $\Z$-grading preserving morphisms $\chi_{+,k}^{(i,j)}$ from $\c\acute{}(\Gamma^R_k)$ to $\c\acute{}(\Gamma^R_{k-1})$ and morphisms $\chi_{-,k}^{(i,j)}$ from $\c\acute{}(\Gamma^L_k)$ to $\c\acute{}(\Gamma^L_{k+1})$. Using these morphisms $\chi_{+,k}^{(i,j)}$ and $\chi_{-,l}^{(i,j)}$, we define the chain complex for colored single crossings as a complex of matrix factorizations in the complex category $\kom(\HMF^{gr})$. 
Then we prove that such morphisms $\chi_{+}^{(i,j)}$ and $\chi_{-}^{(i,j)}$ give a homological link invariant. 
That is, we have to prove that the complexes defined by these morphisms $\chi_{+,k}^{(i,j)}$ and $\chi_{-,k}^{(i,j)}$ induce the following isomorphism in the homotopy category of $\kom(\HMF^{gr})$
\begin{eqnarray*}
&&\c\left( \input{r1p-i} \right) \simeq \c\left( \input{r1c-i}\right)\simeq \c\left( \input{r1m-i}\right),\\[-0.1em]
&&\c\left( \input{r2-ijl} \right) \simeq \c\left( \input{r2-ijc}\right)\simeq \c\left( \input{r2-ijr}\right),\\[-0.1em]
&&\c\left( \input{r2-ijlrev} \right) \simeq \c\left( \input{r2-ijcrev}\right),\c\left( \input{r3-ijk} \right) \simeq \c\left( \input{r3-ijkrev}\right).
\end{eqnarray*}
\end{itemize}
\end{str}
\indent The main purpose of this paper is defining a candidate for matrix factorizations in Strategy \ref{stra} (S1).
For the essential colored planar diagrams (Figure $5$) with additional data which is a sequence naturally induced by coloring, we define matrix factorizations, and then we define matrix factorizations for planar diagrams $\Gamma^R_k$ and $\Gamma^L_k$ of Figure $6$ as tensor product of matrix factorizations for essential planar diagrams. 
In the paper \cite{Yone1}, the author defined matrix factorizations for the planar diagrams \input{figdline2}, \input{fig3valent-in} and \input{fig3valent-out} to reconstruct the matrix factorization for \input{figsmoothing2sln}. Since the planar diagrams \input{figdline2}, \input{fig3valent-in} and \input{fig3valent-out} are the same to colored planar diagrams \input{figline2}, \input{fig3valent-2-1-1} and \input{fig3valent-1-1-2}, we can define matrix factorizations for diagrams of Figure $5$ as a generalization of these matrix factorizations.\\
\indent In the notion of the category $\HMF^{gr}$, we show that some matrix factorizations derived from tensor product of the essential matrix factorizations have isomorphisms in $\HMF^{gr}$ corresponding to MOY relations.\\
%%%%%%%%%%%%%%%%%%%%%%%%%%%%%%%%%%%%%%%
\indent
The present paper is organized as follows.\\
In Section $2$, we give basic definitions, properties and theorems related to category of a matrix factorization. However, many things are already known. New results are Theorem \ref{exclude}, which is a generalization to multivariable of Theorem 2 given by Khovanov and Rozansky \cite{KR3}\cite{Yone1}, Corollary \ref{induce-sq1}, Corollary \ref{cor2-10} and Corollary \ref{cor2-11}. 
In Section $3$, we notice facts of $\Z$-graded algebra. After that, we define matrix factorizations for essential colored planar diagrams of Figure $5$ and show that some matrix factorizations for planar diagrams have equivalences corresponding to MOY relation. \\\\
%acknowledge 
{\bf Acknowledgement.} 
The author is grateful to O. Iyama, H. Ochiai and K. Yoshida, Y. Yoshino, especially A. Tsuchiya and L. Rozansky for their helpful comments.
This work is partly supported by the Grant-in-Aid for JSPS Fellows (20-2330) from Japan Society for the Promotion of Science.

%%%%%%%%%%%%%%%%%%%%%%%%%%%%%%%%%%%%%%%%%%%%%%%%%%%%%%%%%%%%%%%%%%%%%%%%%%%%%%%%%%%%%%%%%%%%%%%%%%%%%%%%
%
%
% section 2 : Category of matrix factorization
% subsection 2.1 : Matrix factorization
%
%
%%%%%%%%%%%%%%%%%%%%%%%%%%%%%%%%%%%%%%%%%%%%%%%%%%%%%%%%%%%%%%%%%%%%%%%%%%%%%%%%%%%%%%%%%%%%%%%%%%%%%%%%

\section{Category of $\Z$-graded matrix factorization}
\subsection{$\Z$-graded matrix factorization}
\indent
Let $R$ be a $\Z$-graded polynomial ring over $\Q$, let $M_0$,$M_1$ be free $\Z$-graded $R$-modules permitted infinite rank and let be a homogeneous element $\omega \in R$. 
A {\bf matrix factorization} 
with a potential $\omega \in R$ is a $2$-cyclic complex of $M_0$ and $M_1$,
$$
\overline{M}=\Big( \xymatrix{
*^{M_0}\ar[rrr]^{d_{M_0}}&&&*^{M_1}\ar[rrr]^{d_{M_1}}&&&*^{M_0}
}
\Big)
$$
with the properties $d_{M_1} d_{M_0} = \omega\,\id_{M_0}$ and $d_{M_0} d_{M_1} = \omega\,\id_{M_1}$. We often simply denote this $2$-cyclic complex by $\overline{M}=(M_0,M_1,d_{M_0},d_{M_1})$.\\
\indent
For a $\Z$-graded polynomial ring $R$ over $\Q$ and $\omega \in R$, let $\MF^{gr}_{R,\omega}$ be a {\bf category of a $\Z$-graded matrix factorization} 
whose object is a matrix factorization $\overline{M}=(M_0,M_1,d_{M_0},d_{M_1})$ with the potential $\omega$ 
and whose morphism $\overline{f}=(f_0,f_1)$ from $\overline{M}=(M_0,M_1,d_{M_0},d_{M_1})$ to $\overline{N}=(N_0,N_1,d_{N_0},d_{N_1})$ consists of a pair of $\Z$-grading preserving $R$-module morphisms 
$f_0:M_0 \to N_0$ and $f_1:M_1 \to N_1$ such that $d_{N_0} f_0 = f_1 d_{M_0}$ and 
$d_{N_1} f_1 = f_0 d_{M_1}$. A matrix factorization $\overline{M}\in \ob (\MF^{gr}_{R,\omega})$ is {\bf finite} if $M_0$ and $M_1$ are finite rank $R$-modules.\\
\indent
Let $R$ and $R\acute{}$ be $\Z$-graded polynomial rings over $\Q$ and $S$ be the maximal ring with inclusions from the ring $S$ to $R$ and $R\acute{}$.
We will take the {\bf tensor product} of $R$ and $R\acute{}$ over the maximal ring $S$ always,
\begin{equation*}
R\ostimes R\acute{}=R\oqtimes R\acute{}/\{rs\oqtimes r\acute{}-r\oqtimes sr\acute{}|r\in R,r\acute{}\in R\acute{},s\in S\}.
\end{equation*} 
A $\Z$-grading of the tensor product $R\ostimes R\acute{}$ is naturally induced by one of the tensor product $R\oqtimes R\acute{}$. 

For an $R$-module $M$ and an $R\acute{}$-module $M\acute{}$, these tensor product over $S$ is also defined by
\begin{equation*}
M\ostimes N=M\oqtimes N/\{ms\oqtimes n-m\oqtimes sn|m\in M,n\in N,s\in S\}.
\end{equation*}
Indeed, a $\Z$-grading of the tensor product $M\ostimes N$ is also induced by one of $M\oqtimes N$.\\
\indent
For $\overline{M}=(M_0,M_1,d_{M_0},d_{M_1}) \in \ob(\MF^{gr}_{R,\omega})$ and $\overline{N}=(N_0,N_1,d_{N_0},d_{N_1}) \in \ob(\MF^{gr}_{R\acute{},\omega\acute{}})$, 
we define the {\bf tensor product} of matrix factorizations $\overline{M}\boxtimes \overline{N} \in \ob(\MF^{gr}_{R\ostimes R\acute{},\omega + \omega \acute{}})$ by
\begin{eqnarray*}
\overline{M}\boxtimes \overline{N} &:=&
\left(\begin{array}{c}
		M_0\ostimes N_0\\
		M_1\ostimes N_1
	\end{array},
	\begin{array}{c}
		M_1\ostimes N_0\\
		M_0\ostimes N_1
	\end{array},
	\left(
		\begin{array}{cc}
			d_{M_0}&-d_{N_1}\\
			d_{N_0}&d_{M_1}
		\end{array}
	    \right),
	    \left(
		\begin{array}{cc}
			d_{M_1}&d_{N_1}\\
			-d_{N_0}&d_{M_0}
		\end{array}
	    \right)\right)
,\\
&=&\left(
\xymatrix{
*^{\left(
	\begin{array}{c}
		M_0\ostimes N_0\\
		M_1\ostimes N_1
	\end{array}
    \right)}
\ar[rrr]_{\left(
		\begin{array}{cc}
			d_{M_0}&-d_{N_1}\\
			d_{N_0}&d_{M_1}
		\end{array}
	    \right)}&&&
*^{\left(
	\begin{array}{c}
		M_1\ostimes N_0\\
		M_0\ostimes N_1
	\end{array}
    \right)}
\ar[rrr]_{\left(
		\begin{array}{cc}
			d_{M_1}&d_{N_1}\\
			-d_{N_0}&d_{M_0}
		\end{array}
	    \right)}
&&&*^{\left(
	\begin{array}{c}
		M_0\ostimes N_0\\
		M_1\ostimes N_1
	\end{array}
    \right)}
}
\right)
\end{eqnarray*}
where 
\begin{center}
$\left(
		\begin{array}{cc}
			d_{M_0}&-d_{N_1}\\
			d_{N_0}&d_{M_1}
		\end{array}
	    \right)$ 
and
$\left(
		\begin{array}{cc}
			d_{M_1}&d_{N_1}\\
			-d_{N_0}&d_{M_0}
		\end{array}
	    \right) $
\end{center}
simply denote 
\begin{center}
$\left(
		\begin{array}{cc}
			d_{M_0} \ostimes \id_{N_0}&- \id_{M_1}\ostimes d_{N_1}\\
			\id_{M_0}\ostimes d_{N_0}&d_{M_1}\ostimes \id_{N_1}
		\end{array}
	    \right)$ 
and
$\left(
		\begin{array}{cc}
			d_{M_1}\ostimes \id_{N_0}&\id_{M_0}\ostimes d_{N_1}\\
			-\id_{M_1}\ostimes d_{N_0}&d_{M_0}\ostimes \id_{N_1}
		\end{array}
	    \right) $.
\end{center}

\indent
This tensor product $\boxtimes$ is commutative and associative. Moreover, there is the unit object for the tensor product.\\

\begin{pro}\label{com-ass}%proposition
{\rm (1)}For $\overline{M} \in \ob(\MF^{gr}_{R,\omega})$ and $\overline{N} \in \ob(\MF^{gr}_{R\acute{},\omega\acute{}})$, 
there is an isomorphism in $\MF^{gr}_{R\ostimes R\acute{},\omega + \omega \acute{}}$
$$\overline{M}\boxtimes \overline{N} \simeq \overline{N}\boxtimes \overline{M}.$$ 
{\rm (2)}For $\overline{L} \in \ob(\MF^{gr}_{R,\omega})$, $\overline{M} \in \ob(\MF^{gr}_{R\acute{},\omega\acute{}})$ 
and $\overline{N} \in \ob(\MF^{gr}_{R \acute{}\,\acute{},\omega \acute{}\,\acute{}})$, 
there is an isomorphism in $\MF^{gr}_{R\ostimes R\acute{}\osatimes R \acute{}\,\acute{},\omega + \omega\acute{} +\omega \acute{}\,\acute{}}$ 
$$(\overline{L}\boxtimes \overline{M})\boxtimes \overline{N} \simeq \overline{L}\boxtimes (\overline{M}\boxtimes \overline{N}).$$ 
\end{pro}
\begin{proof}
See \cite{Yone1}.
\end{proof}
\begin{rem}
As from here, $\overline{M_1}\boxtimes\overline{M_2}\boxtimes\overline{M_3}\boxtimes\cdots\boxtimes\overline{M_n}$ is expressed as
\begin{equation*}
\overline{M_1}\boxtimes\overline{M_2}\boxtimes\overline{M_3}\boxtimes\cdots\boxtimes\overline{M_n}
=(\cdots((\overline{M_1}\boxtimes\overline{M_2})\boxtimes\overline{M_3})\boxtimes\cdots)\boxtimes\overline{M_n}.
\end{equation*}
\end{rem}
\begin{pro}\label{identity}%proposition
The matrix factorization $(\xymatrix{R \ar[r]^0&0\ar[r]^0&R})$ is the unit object for tensor product to any object in $\MF^{gr}_{R,\omega}$. In brief, 
for any matrix factorization $\overline{M} \in \ob(\MF^{gr}_{R,\omega})$ we have
$$
\overline{M} \boxtimes (\xymatrix{R \ar[r]^0&0\ar[r]^0&R}) \simeq \overline{M}.
$$
\end{pro}

\begin{proof}
See \cite{Yone1}.
\end{proof}

$\overline{R}$ often simply denotes 
$$
\overline{R}:=\xymatrix{(R\ar[r]^0&0\ar[r]^0&R)}.
$$

\indent
The {\bf translation functor} $\left< 1 \right>$ changes a matrix factorization $\overline{M} = (M_{0},M_{1},d_{M_0},d_{M_1})\in \MF^{gr}_{R,\omega}$ into 
\begin{eqnarray*}
\overline{M} \left< 1 \right> &=&
(
M_{1},M_{0},-d_{M_1},-d_{M_0})\in \MF^{gr}_{R,\omega}\\
\Big(&=&\left(
\xymatrix{
*{M_{1}}\ar[rr]^{-d_{M_1}}&&*{M_{0}}\ar[rr]^{-d_{M_0}}&&*{M_{1}}
}
\right)\Big)
.
\end{eqnarray*}
The functor $\left< 2 \right> (= \left< 1 \right>^2)$ is the identity functor.
\begin{pro}%proposition
For $\overline{M} \in \ob(\MF^{gr}_{R,\omega})$ and $\overline{N} \in \ob(\MF^{gr}_{R\acute{},\omega\acute{}})$, 
there is an isomorphism in $\MF^{gr}_{R\ostimes R\acute{},\omega + \omega \acute{}}$
\begin{eqnarray*}
(\overline{M}\boxtimes \overline{N})\left< 1 \right> &=& (\overline{M}\left< 1 \right> )\boxtimes \overline{N}\\[-0.1em]
&\simeq& \overline{M}\boxtimes (\overline{N}\left< 1 \right> ).
\end{eqnarray*}
\end{pro}
\begin{proof}
See \cite{Yone1}.
\end{proof}
\indent
The morphism $\overline{f}=(f_1,f_2): \overline{M} \to \overline{N} \in \mor (\MF^{gr}_{R,\omega})$ is {\bf null-homotopic} 
if morphisms $h_{0}:M_0 \to N_1$ and $h_{1}:M_1 \to N_0$ exist 
such that $f_0 = h_1 d_{M_0} + d_{N_1} h_0$ and $f_1 = h_0 d_{M_1} + d_{N_0} h_1$.
Two morphisms $\overline{f},\overline{g}: \overline{M} \to \overline{N} \in \mor (\MF^{gr}_{R,\omega})$ are {\bf homotopic} 
if $\overline{f} - \overline{g}$ is null-homotopic. In this case, we denote $\overline{f} \sim \overline{g}$. Two matrix factorization $\overline{M}$ and $\overline{N}\in\ob (\MF^{gr}_{R,\omega})$ are {\bf homotopy equivalence} if there are $\Z$-grading preserving morphisms $f:\overline{M}\to \overline{N}$ and $g:\overline{N}\to \overline{M}$ such that $fg\sim \id_{\overline{N}}$ and $gf\sim \id_{\overline{M}}$.\\
\indent
Let $\HMF^{gr}_{R,\omega}$ be the quotient category of $\MF^{gr}_{R,\omega}$ which has the same objects to $\MF^{gr}_{R,\omega}$
and has morphisms of $\mor(\MF^{gr}_{R,\omega})$ modulo null-homotopic. The category $\HMF^{gr}_{R,\omega}$ is called the {\bf homotopy category} of $\MF^{gr}_{R,\omega}$. It is obvious that a homotopy equivalence in $\MF^{gr}_{R,\omega}$ is an isomorphism in $\HMF^{gr}_{R,\omega}$.

\indent A matrix factorization in $\MF^{gr}_{R,\omega}$ is called {\bf contractible} 
if it is isomorphic in $\HMF^{gr}_{R,\omega}$ to the zero matrix factorization
$$
\Big(
\xymatrix{
*{0}\ar[rr]&&*{0}\ar[rr]&&*{0} 
}
\Big).
$$
\indent
A {\bf $\Z$-grading shift} $\{ n \}$ ($n \in \Z$) is an operator up $n$-grading.
That is, for a $\Z$-graded $R$-module $M$ with the $\Z$-graded decomposition $\oplus M^i$ ($M^i$:$\Q$ vector space with $i$-grading), $M\{n\}$ is defined by
$$
M\{n\}=\oplus M^{i+n}.
$$
The $\Z$-grading shift $\{ n \}$ turns the matrix factorization $\overline{M}=(M_0,M_1,d_{M_0},d_{M_1})$ into 
\begin{eqnarray*}
\overline{M}\{n\}&=&(M_0\{ n \},M_1\{ n\},d_{M_0},d_{M_1})\\
&=&\Big( \xymatrix{ M_0\{ n \} \ar[r]^{ d_{M_0} }   &   M_1\{ n\}   \ar[r]^{ d_{M_1} }   &   M_0\{ n\} }\Big).
\end{eqnarray*}
\begin{pro}\label{functor1}%proposition
For $\overline{M} \in \ob(\MF^{gr}_{R,\omega})$ and $\overline{N} \in \ob(\MF^{gr}_{R\acute{},\omega\acute{}})$, 
there is an equality in $\MF^{gr}_{R\ostimes R\acute{},\omega + \omega \acute{}}$
\begin{eqnarray*}
(\overline{M}\boxtimes \overline{N})\{ n\} &=& (\overline{M}\{ n\} )\boxtimes \overline{N}\\[-0.1em]
&=& \overline{M}\boxtimes (\overline{N}\{ n\} )
\end{eqnarray*}
\end{pro}
\begin{proof}
We find that these objects are really identical by definition. 
\end{proof}
If the potential $\omega$ equals $0$ then a matrix factorization $\overline{M} \in\ob(\MF^{gr}_{R,0})$ satisfies the boundary condition $d_{M_0} d_{M_1} = d_{M_1} d_{M_0} =0 $.  
Therefore, there is homology functor $\H$ from $\MF^{gr}_{R,0}$ to a category of a $\Z \oplus \Z_2$-graded homology group as follows,
$$ \H(\overline{M})=\bigoplus_{j \in \Z ,k\in \Z_2 } \H^{j,k}(\overline{M}),$$
where $k$ is a complex grading of the matrix factorization $\overline{M}$ (i.e. $k=0$ or $1$) 
and $j$ is a $\Z$-grading derived from the $\Z$-graded modules of the matrix factorization $\overline{M}$.
The {\bf Euler characteristic} $\overline{\chi}$ of $H(\overline{M})$ is defined by
\begin{equation*}
\overline{\chi}(\H(\overline{M})) = \sum_{ j \in \Z, k \in \Z_2} \dim_{\Q} \H^{j,k}(\overline{M})  q^{j} \label{euler}.
\end{equation*}

%%%%%%%%%%%%%%%%%%%%%%%%%%%%%%%%%%%%%%%%%%%%%%%%%%%%%%%%%%%%%%%%%%%%%%%%%%%%%%%%%%%%%%%%%%%%%%%%%%%%%%%%
%
%
% subsection 2.4 : Koszul matrix factorization
%
%
%%%%%%%%%%%%%%%%%%%%%%%%%%%%%%%%%%%%%%%%%%%%%%%%%%%%%%%%%%%%%%%%%%%%%%%%%%%%%%%%%%%%%%%%%%%%%%%%%%%%%%%%

\subsection{Koszul matrix factorization}
Let $R$ be a $\Z$-graded polynomial ring over $\Q$. 
For homogeneous $\Z$-grading polynomials $a$, $b \in R$ and a $\Z$-graded $R$-module $M$, 
we define a matrix factorization $ K(a;b)_{M} $
with the potential $a b$ by
\begin{eqnarray*}
K(a;b)_{M}&=& 
(M,M\{ \frac{1}{2}(\, \deg (b)-\deg (a)\, )\},a,b)\\
&=& 
\Big(
\xymatrix{
*{M}\ar[rr]^(.3){a}&&*{M\{ \frac{1}{2}(\, \deg (b)-\deg (a)\, )\}}\ar[rr]^(.7){b}&&*{M}
}
\Big)
,
\end{eqnarray*}
where $\deg (a)$ and $\deg (b)$ are $\Z$-gradings of the polynomials $a$, $b \in R$.
In general, for sequences $\mathbf{a}={}^t(a_1, a_2, \ldots, a_k)$, $\mathbf{b}={}^t(b_1, b_2, \ldots, b_k)$ of homogeneous polynomials in $R$ and $R$-module $M$, 
a matrix factorization $K\left( \mathbf{a} ; \mathbf{b} \right)_{M}$ with the potential $\sum_{i=1}^k a_i b_i$ is defined by 
\[
K\left( \mathbf{a} ; \mathbf{b} \right)_{M}
=
\mathop{\boxtimes}_{i=1}^k K(a_i;b_i)_{R}\boxtimes(M,0,0,0)
.
\]
This matrix factorization is called a {\bf Koszul matrix factorization} \cite{KR1}.
\begin{rem}
Let $R$ be a $\Z$-graded polynomial ring over $\Q$ and let $R_y$ be the $\Z$-graded polynomial ring $R[y]$. For polynomials $a$ and $b$ in $R$,
$K(a;b)_{R_y}$ is a matrix factorization of $R_y$-modules with rank $1$ in $\MF^{gr}_{R_y,ab}$ 
and a matrix factorization of infinite rank $R$-modules in $\MF^{gr}_{R,ab}$ besides.
\end{rem}
\begin{pro}\label{equiv2}%proposition
Let $c$ be a non-zero element in $\Q$.
There is the following isomorphism in $\MF^{gr}_{R,a b}$
$$
K(a;b)_M \simeq K(c a;c^{-1} b)_M. 
$$
\end{pro}
\begin{proof}
See \cite{Yone1}.
\end{proof}
\begin{pro}\label{functor2}%proposition
\begin{eqnarray*}
K(a;b)_{M}\,\left< 1 \right> &=& K(-b;-a)_{M}\,\{ \frac{1}{2}(\, \deg (b) -\deg (a)\, ) \} \\[-0.1em]
&\simeq&K(b;a)_{M}\,\{ \frac{1}{2}(\, \deg (b) -\deg (a)\, ) \} 
\end{eqnarray*}
\end{pro}
\begin{proof}
See \cite{Yone1}.
%The first equation is obvious by definition.
%It is easy to find the second equivalence using Proposition \ref{equiv2} as $c=-1$.
%%By definition, we have 
%%\begin{eqnarray*}
%%K(a;b)_{R}\left< 1 \right> &=& \Big(\xymatrix{R \ar[rr]^(.3){a}&& R\{ \frac{1}{2}(\deg b -\deg a) \} \ar[rr]^(.7){b}&& R}\Big)\,\left< 1 \right> \\
%%&=& \Big(\xymatrix{R\{ \frac{1}{2}(\deg b -\deg a) \} \ar[rr]^(.7){-b}&& R \ar[rr]^(.3){-a}&& R\{ \frac{1}{2}(\deg b -\deg a) \} }\Big)\\
%%&=& \Big(\xymatrix{R \ar[rr]^(.3){-b}&& R\{ \frac{1}{2}(\deg a -\deg b) \} \ar[rr]^(.7){-a}&& R}\Big)\,\{ \frac{1}{2}(\deg b -\deg a) \} \\
%%&=& K(-b;-a)_{R}\,\{ \frac{1}{2}(\deg b -\deg a) \} .
%%\end{eqnarray*}
\end{proof}
%\begin{pro}
%Let $a$ and $b$ be polynomials with a homogeneous $\Z$-grading in $R$ and let $M$ be a finite $\Z$-graded $R$-module. There is the following isomorphism in $\MF^{gr}_{R,-ab}$,
%\begin{eqnarray*}
%K(a;b)_{M \bullet}\simeq K(b;-a)_{M^{\ast}}.
%\end{eqnarray*}
%\end{pro}
%\begin{proof}
%We have
%\begin{eqnarray*}
%K(a;b)_{M \bullet}&=&\overline{\Hom_R}(K(a;b)_{M},\overline{R})\\
%&\simeq&(\xymatrix{\Hom_R(M,R)\ar[r]^(0.325){b}&\Hom_R(M\{\frac{1}{2}(\deg b-\deg a)\},R)\ar[r]^(0.65){-a}&\Hom_R(M,R)}).
%\end{eqnarray*}
%$\Hom_R(M,R)$ and $\Hom_R(M\{\frac{1}{2}(\deg b-\deg a)\},R)$ are isomorphic to $M^{\ast}$ and $M^{\ast}\{\frac{1}{2}(\deg a-\deg b)\}$ respectively.
%Thus we have
%\begin{eqnarray*}
%K(a;b)_{M \bullet}&\simeq&(\xymatrix{M^{\ast}\ar[r]^(0.25){b}&M^{\ast}\{\frac{1}{2}(\deg a-\deg b)\}\ar[r]^(0.7){-a}&M^{\ast}})\\
%&\simeq&K(b;-a)_{M^{\ast}}.
%\end{eqnarray*}
%\end{proof}

\begin{pro}\label{equiv}%proposition
Let $a_i$ and $b_i$ be homogeneous $\Z$-grading polynomials such that $\deg (a_1) +\deg (b_1) =\deg (a_2) +\deg (b_2)$ and let $\lambda_i$ {\rm ($i=1,2$)} be homogeneous $\Z$-grading polynomials such that $\deg ( \lambda_1 )=\deg (a_2) - \deg (a_1)$ and $\deg ( \lambda_2 )=-\deg (b_1) + \deg (a_2)$. \\
{\rm (1)} There is the following isomorphism in $\MF^{gr}_{R,a_1 b_1 + a_2 b_2}$ 
$$
K\left(\left(
\begin{array}{c}
	 a_1\\
	 a_2
\end{array}
\right);
\left(
\begin{array}{c}
	 b_1\\
	 b_2
\end{array}
\right)\right)_{M}
\simeq
K\left(\left(
\begin{array}{c}
	 a_1\\
	 a_2 + \lambda_1 a_1
\end{array}
\right);
\left(
\begin{array}{c}
	 b_1 - \lambda_1 b_2\\
	 b_2
\end{array}
\right)\right)_{M}.
$$
{\rm (2)} There is the following isomorphism in $\MF^{gr}_{R,a_1 b_1 + a_2 b_2}$ 
$$
K\left(\left(
\begin{array}{c}
	 a_1\\
	 a_2
\end{array}
\right);
\left(
\begin{array}{c}
	 b_1\\
	 b_2
\end{array}
\right)\right)_{M}
\simeq
K\left(\left(
\begin{array}{c}
	 a_1 - \lambda_2 b_2\\
	 a_2 + \lambda_2 b_1
\end{array}
\right);
\left(
\begin{array}{c}
	 b_1\\
	 b_2
\end{array}
\right)\right)_{M}.
$$

\end{pro}

\begin{proof}
See \cite{Ras}\cite{Wu}.
\end{proof}
\begin{thm}\label{reg-eq}%theorem
Let $R$ be a $\Z$-graded polynomial ring and let $a_i$, $b_i$ and $b_i\acute{}$ ($i=1,\cdots ,m$) be polynomials in $R$.
If $a_1,\cdots ,a_m$ $\in R$ form a regular sequence and
\begin{equation*}
\sum_{i=1}^m a_ib_i=\sum_{i=1}^m a_ib_i\acute{} \,(=:\omega),
\end{equation*}
there is the following isomorphism in $\MF^{gr}_{R,\omega}$
\begin{equation*}
k\left(
\left(
\begin{array}{c}
a_1\\
\vdots\\
a_m
\end{array}
\right);
\left(
\begin{array}{c}
b_1\\
\vdots\\
b_m
\end{array}
\right)
\right)_M
\simeq
k\left(
\left(
\begin{array}{c}
a_1\\
\vdots\\
a_m
\end{array}
\right);
\left(
\begin{array}{c}
b_1\acute{}\\
\vdots\\
b_m\acute{}
\end{array}
\right)
\right)_M.
\end{equation*}
\end{thm}
\begin{proof}
See Theorem 2.1 on \cite{KR3}
\end{proof}

\begin{cor}\label{induce-sq1}%Corollary
Put $R=\Q[x_1,x_2,\cdots ,x_k]$ and $R_y=R[y]\left/\left<y^l+\alpha_1y^{l-1} +\alpha_2y^{l-2} +\cdots +\alpha_l\right>\right.$ ($\alpha_i\in R$).\\
($1$)Let $a_i$ be a polynomial $\in R_y$($i=1,\cdots ,m$), $b_i$ be a polynomial $\in R$ ($i=2,\cdots ,m$) and let $b_1$, $\beta$ be polynomials $\in R_y$ with the property $(y+\beta)b_1\in R$.
If these polynomials hold the following conditions;
\begin{itemize}
\item[(i)]$(y+\beta)b_1$, $b_2$, $\cdots$, $b_m$ form a regular sequence in $R$,
\item[(ii)] $\displaystyle a_1b_1(y+\beta)+\sum_{i=2}^m a_ib_i$ $(=:\omega)$ $\in R$,
\end{itemize}
then there are polynomials $a_i\acute{} \in R$ ($i=1,\cdots ,m$) to give the following isomorphism in $\MF_{R,\omega}$
\begin{equation*}
K\left(\left(
\begin{array}{c}
(y+\beta)a_1\\
a_2\\
\vdots\\
a_m
\end{array}
\right);
\left(
\begin{array}{c}
b_1\\
b_2\\
\vdots\\
b_m
\end{array}
\right)\right)_{R_y}
\simeq
K\left(\left(
\begin{array}{c}
(y+\beta)a_1\acute{}\\
a_2\acute{}\\
\vdots\\
a_m\acute{}
\end{array}
\right);
\left(
\begin{array}{c}
b_1\\
b_2\\
\vdots\\
b_m
\end{array}
\right)\right)_{R_y}.
\end{equation*}
\\
\noindent
($2$)Let $a_i$ be a polynomial $\in R_y$ ($i=1,\cdots ,m$), $b_i$ be a polynomial $\in R$ ($i=1,\cdots ,m$) and $\beta$ be a polynomial $\in R$.
If these polynomials hold the following conditions;
\begin{itemize}
\item[(i)]$b_1$, $b_2$, $\cdots$, $b_m$ form a regular sequence in $R$,
\item[(ii)] $\displaystyle a_1b_1(y+\beta)+\sum_{i=2}^m a_ib_i$ $(=:\omega\acute{})$ $\in R$,
\end{itemize}
then there are polynomials $a_1\acute{} \in R_y$ and $a_i\acute{} \in R$ ($i=2,\cdots ,m$) to give the following isomorphism in $\MF_{R,\omega\acute{}}$
\begin{equation*}
K\left(\left(
\begin{array}{c}
a_1\\
a_2\\
\vdots\\
a_m
\end{array}
\right);
\left(
\begin{array}{c}
b_1(y+\beta)\\
b_2\\
\vdots\\
b_m
\end{array}
\right)\right)_{R_y}
\simeq
K\left(\left(
\begin{array}{c}
a_1\acute{}\\
a_2\acute{}\\
\vdots\\
a_m\acute{}
\end{array}
\right);
\left(
\begin{array}{c}
b_1(y+\beta)\\
b_2\\
\vdots\\
b_m
\end{array}
\right)\right)_{R_y}.
\end{equation*}
\end{cor}

\begin{proof}
($1$)We consider the matrix factorization
\begin{equation*}
K\left(\left(
\begin{array}{c}
a_1\\
a_2\\
\vdots\\
a_m
\end{array}
\right);
\left(
\begin{array}{c}
(y+\beta)b_1\\
b_2\\
\vdots\\
b_m
\end{array}
\right)\right)_{R_y}.
\end{equation*}
Since the sequence $((y+\beta)b_1,b_2,\cdots,b_m)$ is regular, by Corollary \ref{reg-eq} there are polynomials $a_i\acute{} \in R$ ($i=1,\cdots ,m$) to give the following matrix factorization,
\begin{equation*}
K\left(\left(
\begin{array}{c}
a_1\acute{}\\
a_2\acute{}\\
\vdots\\
a_m\acute{}
\end{array}
\right);
\left(
\begin{array}{c}
(y+\beta)b_1\\
b_2\\
\vdots\\
b_m
\end{array}
\right)\right)_{R_y}.
\end{equation*}
It is obvious that the polynomial degree of $b_1$ is less than one of $(y+\beta)b_1$. Then we have the isomorphism
\begin{equation*}
K\left(\left(
\begin{array}{c}
(y+\beta)a_1\\
a_2\\
\vdots\\
a_m
\end{array}
\right);
\left(
\begin{array}{c}
b_1\\
b_2\\
\vdots\\
b_m
\end{array}
\right)\right)_{R_y}
\simeq
K\left(\left(
\begin{array}{c}
(y+\beta)a_1\acute{}\\
a_2\acute{}\\
\vdots\\
a_m\acute{}
\end{array}
\right);
\left(
\begin{array}{c}
b_1\\
b_2\\
\vdots\\
b_m
\end{array}
\right)\right)_{R_y}.
\end{equation*}
\\
\noindent($2$)The polynomial $a_1$ can be described as 
\begin{equation}
\label{quoti-poly}a_1=\gamma_{0}y^{l-1}+\gamma_{1}y^{l-2}+\cdots +\gamma_{l-1} \hspace{1cm}(\gamma_i \in R).
\end{equation}
First we consider the matrix factorization
\begin{equation*}
K\left(\left(
\begin{array}{c}
a_1(y+\beta)\\
a_2\\
\vdots\\
a_m
\end{array}
\right);
\left(
\begin{array}{c}
b_1\\
b_2\\
\vdots\\
b_m
\end{array}
\right)\right)_{R_y}.
\end{equation*}
By assumption, the second sequence $(b_1,b_2,\cdots,b_m)$ is regular and potential of the matrix factorization is in $R$. 
Using Theorem \ref{reg-eq} we find that there are polynomials $c_1$ and $a_i\acute{}$ ($i=2,\cdots ,m$) in $R$ to give an isomorphism between the above matrix factorization and the following matrix factorization
\begin{equation*}
K\left(\left(
\begin{array}{c}
c_1\\
a_2\acute{}\\
\vdots\\
a_m\acute{}
\end{array}
\right);
\left(
\begin{array}{c}
b_1\\
b_2\\
\vdots\\
b_m
\end{array}
\right)\right)_{R_y}.
\end{equation*}
Therefore, $a_1(y+\beta)-c_1$ can be described as an $R_y$ coefficient linear sum of $b_2$, $\cdots$, $b_{m-1}$, $b_m$. Since, by Equation (\ref{quoti-poly}),
\begin{eqnarray*}
a_1(y+\beta)&=&\gamma_{0}y^{l}+(\gamma_{1}+\beta\gamma_{0})y^{l-1}+(\gamma_{2}+\beta\gamma_{1})y^{l-2}+\cdots +(\gamma_{l-1}+\beta\gamma_{l-2})y+\beta\gamma_{l-1}\\[-0.1em]
&=&(\gamma_{1}+\beta\gamma_{0}-\alpha_1\gamma_{0})y^{l-1}+\cdots +(\gamma_{l-1}+\beta\gamma_{l-2}-\alpha_{l-1}\gamma_{0})y+\beta\gamma_{l-1}-\alpha_l\gamma_{0}
\end{eqnarray*}
and $b_2$, $\cdots$, $b_{m-1}$, $b_{m}$ form a regular sequence, the polynomials $\gamma_{i}+\beta\gamma_{i-1}-\alpha_{i}\gamma_{0}$ ($i=1,\cdots ,l-1$) are described as an $R$ coefficient linear sum of $b_2$, $\cdots$, $b_{m-1}$, $b_m$,
\begin{equation}\label{gamma-eq}
\gamma_{i}+\beta\gamma_{i-1}-\alpha_{i}\gamma_{0}=\sum_{j=2}^{m}s_{i,j}b_j \hspace{1cm}(i=1,\cdots ,l-1),
\end{equation}
and $\beta\gamma_{l-1}-\alpha_l\gamma_{0}$ is described as
\begin{equation*}
\beta\gamma_{l-1}-\alpha_l\gamma_{0}=c_1+\sum_{j=2}^{m}s_{l,j}b_j.
\end{equation*}
Moreover equations of (\ref{gamma-eq}) change into the following equations by linear transform,
\begin{equation*}
\gamma_i=(\sum_{j=0}^{i}(-1)^{i-j}\alpha_j\beta^{i-j})\gamma_0+\sum_{j=2}^{m}\widetilde{s_{i,j}}b_j\hspace{1cm}(i=1,\cdots ,l-1),
\end{equation*}
where $\alpha_0=1$ and $\displaystyle\widetilde{s_{i,j}}= \sum_{k=1}^{i}(-1)^k\beta^k s_{i-k,j}$.\\
Thus we have the equivalence between matrix factorizations as follows,
\begin{equation*}
K\left(\left(
\begin{array}{c}
a_1\\
a_2\\
\vdots\\
a_m
\end{array}
\right);
\left(
\begin{array}{c}
b_1(y+\beta)\\
b_2\\
\vdots\\
b_m
\end{array}
\right)\right)_{R_y}
\simeq
K\left(\left(
\begin{array}{c}
a_1\acute{}\\
a_2-\sum_{i=1}^{l-1}\widetilde{s_{i,2}}b_1y^{l-1-i}\\
\vdots\\
a_m-\sum_{i=1}^{l-1}\widetilde{s_{i,m}}b_1y^{l-1-i}
\end{array}
\right);
\left(
\begin{array}{c}
b_1(y+\beta)\\
b_2\\
\vdots\\
b_m
\end{array}
\right)\right)_{R_y},
\end{equation*}
where
\begin{equation*}
a_1\acute{}=\sum_{i=0}^{l-1}\left( \sum_{j=0}^{i}(-1)^{i-j}\alpha_j\beta^{i-j}\right)\gamma_0y^{l-1-i}.
\end{equation*}
It is obvious to find the polynomial $a_1\acute{}\,(y+\beta)$ is in $R$. Since the polynomials $b_2$, $\cdots$, $b_{m-1}$, $b_m$ form a regular sequence and the polynomial $\sum_{j=2}^{m}(a_j-\sum_{i=1}^{l-1}\widetilde{s_{i,j}}b_1y^{l-1-i})b_j=\omega\acute{}-a_1\acute{}\,(y+\beta)$ is in $R$,
polynomials $a_i\acute{} \in R$ ($i=2,\cdots ,m$) exist to give the following isomorphism 
\begin{equation*}
K\left(\left(
\begin{array}{c}
a_1\acute{}\\
a_2-\sum_{i=1}^{l-1}\widetilde{s_{i,2}}b_1y^{l-1-i}\\
\vdots\\
a_m-\sum_{i=1}^{l-1}\widetilde{s_{i,m}}b_1y^{l-1-i}
\end{array}
\right);
\left(
\begin{array}{c}
b_1(y+\beta)\\
b_2\\
\vdots\\
b_m
\end{array}
\right)\right)_{R_y}
\simeq
K\left(\left(
\begin{array}{c}
a_1\acute{}\\
a_2\acute{}\\
\vdots\\
a_m\acute{}
\end{array}
\right);
\left(
\begin{array}{c}
b_1(y+\beta)\\
b_2\\
\vdots\\
b_m
\end{array}
\right)\right)_{R_y}.
\end{equation*}
\end{proof}
\begin{rem}%Remark
Put $R_y=R[y]\left/\left<y^l+\alpha_1y^{l-1} +\alpha_2y^{l-2} +\cdots +\alpha_l\right>\right.$ ($\alpha_i\in R$). If the variable $y$ remains in a polynomial $p$ then the multiplication map $p$ to $R_y$ is complicated as an $R$-module morphism. However, if the variable $y$ dose not exist in a polynomial $p$ then the multiplication map $p$ is simply a diagonal map as an $R$-module morphism,
\begin{equation*}
\xymatrix{
R_y\ar[r]^{p}&R_y
}
=
\xymatrix{
\hbox{$
\left(
\begin{array}{r}
\alpha_0 R\\
\alpha_1 R\\
\vdots\\
\alpha_{l-1} R
\end{array}
\right)
$}
\ar[rrrr]_{
\hbox{$
\left(
\begin{array}{cccc}
p&0&\cdots&0\\
0&p&&0\\
\vdots&&\ddots&\vdots\\
0&0&\cdots&p
\end{array}
\right)
$}
}&&&&
\hbox{
$\left(
\begin{array}{r}
\alpha_0 R\\
\alpha_1 R\\
\vdots\\
\alpha_{l-1} R
\end{array}
\right),
$}
}
\end{equation*}
where $\alpha_0$, $\alpha_1$, $\cdots$, $\alpha_{l-1}$ form a basis of $R_y$ as an $R$-module.
\end{rem}
In next section, Corollary \ref{induce-sq1} is useful for decomposing a matrix factorization into its direct sum of matrix factorizations.\\

\indent The following theorem is a generalization to multivariable of Theorem 2.2 given by Khovanov and Rozansky \cite{KR3}\cite{Yone1}.

\begin{thm}\label{exclude}%theorem
We put $R=\Q [\underline{x}]$ and $R_{\bf y} = R [\underline{y}]$, where $\underline{x}=(x_1,x_2,\cdots ,x_l)$ and $\underline{y}=(y_1,y_2,\cdots ,y_m)$. 
If $\mathbf{a}={}^t( a_1 , a_2 , \cdots , a_k )$ and $\mathbf{b}={}^t( b_1 , b_2 , \cdots , b_k )$, 
where $a_i$ and $b_i$ are homogeneous polynomials in $R_{\bf y}$, satisfy the following conditions
\begin{itemize}
\item[(i)] $\sum_{i=1}^{k}a_i b_i \,(=:\omega) \in R$,\\
\item[(ii)] There exists $b_j$ which can be described by $c y_1^{i_1}y_2^{i_2}\cdots y_m^{i_m} + p$, where
%\begin{center}
$c$ is a constant and $p$ is a polynomial of $R_{\bf y}$ whose every monomials are not zero in the quotient ring $R_{\bf y}/\left< y_1^{i_1}y_2^{i_2}\cdots y_m^{i_m}\right>$,
%\end{center}
\end{itemize}
then there is the following isomorphism in $\HMF^{gr}_{R,\omega }$,
$$
K(\mathbf{a};\mathbf{b})_{R_{\bf y}} \simeq K(\stackrel{j}{\check{\mathbf{a}}};\stackrel{j}{\check{\mathbf{b}}})_{R_{\bf y}/\left< b_j\right> },
$$
where $\stackrel{j}{\check{\mathbf{a}}}$ and $\stackrel{j}{\check{\mathbf{b}}}$ are the sequences omitted the $i$-th polynomial of $\mathbf{a}$ and $\mathbf{b}$.
\end{thm}

\begin{proof}
We can prove this theorem to repeat proof by M.Khovanov and L.Rozansky in \cite{KR3}\cite{Yone1}.
By Proposition \ref{equiv2} and the above assumption, we can describe the left-hand matrix factorization 
$$
K(\mathbf{a};\mathbf{b})_{R_{\bf y}} \simeq K(\stackrel{i}{\check{\mathbf{a}}};\stackrel{i}{\check{\mathbf{b}}})_{R_{\bf y}} \boxtimes K(a_j;y_1^{i_1}y_2^{i_2}\cdots y_m^{i_m} + p)_{R_{\bf y}}.
$$
The Koszul matrix factorization $K(\stackrel{i}{\check{\mathbf{a}}};\stackrel{i}{\check{\mathbf{b}}})_{R_{\bf y}}$ is described as $(R_{{\bf y}}^r,R_{{\bf y}}^r,D_0,D_1)$, where $D_0$, $D_1 \in \M_r(R_{\bf y})$ ($r=2^{k-2}$).
Then we have
$$
K(\mathbf{a};\mathbf{b})_{R_{\bf y}}=
\left(
\xymatrix{
*{
\begin{array}{c}
R_{\bf y}^r\\
\oplus\\
R_{\bf y}^r
\end{array}
}
\ar[rrr]^{D_0 \acute{}}&&&
*{
\begin{array}{c}
R_{\bf y}^r\\
\oplus\\
R_{\bf y}^r
\end{array}
}
\ar[rrr]^{D_1 \acute{}
}&&&
*{
\begin{array}{c}
R_{\bf y}^r\\
\oplus\\
R_{\bf y}^r
\end{array}
}
}
\right)
,$$
where $D_0 \acute{} =\left(
\begin{array}{cc}
D_0&(-y_1^{i_1}y_2^{i_2}\cdots y_m^{i_m} -p) \id_{R_{\bf y}^r}\\
a_j \id_{R_{\bf y}^r}&D_1
\end{array}
\right)$
,
$D_1 \acute{}=\left(
\begin{array}{cc}
D_1&(y_1^{i_1}y_2^{i_2}\cdots y_m^{i_m}+ p) \id_{R_{\bf y}^r}\\
-a_j \id_{R_{\bf y}^r}&D_0
\end{array}
\right)$.

\indent
The ring $R_{\bf y}$ is split into the direct sum as an $R$-module;
$$
R_{\bf y} \simeq R_{<} \oplus R_{\geq},
$$
where 
$$
\displaystyle R_{<}= \hspace{-1cm}\bigoplus_{\tiny \begin{array}{c}(i_1,i_2,\cdots ,i_m) \in \N_{\geq 0}^m \\ \min \{ i_1-n_1,i_2-n_2,\cdots ,i_m-n_m \} < 0 \end{array}} \hspace{-1cm}R\cdot y_1^{i_1}y_2^{i_2}\cdots y_m^{i_m} \simeq R_{\bf y}/\left<y_1^{n_1}y_2^{n_2}\cdots y_m^{n_m} + p\right>
$$ 
and 
$$
\displaystyle R_{\geq }=\hspace{-0.5cm}\bigoplus_{(i_1,i_2,\cdots ,i_m) \in \N_{\geq 0}^m} \hspace{-0.5cm}R \cdot y_1^{i_1}y_2^{i_2}\cdots y_m^{i_m} (y_1^{n_1}y_2^{n_2}\cdots y_m^{n_m} + p) .
$$\\
The $R$-module morphism $\xymatrix{R_{\bf y} \ar[rrr]^{y_1^{n_1}y_2^{n_2}\cdots y_m^{n_m} + p}&&& R_{\bf y}}$ induces the $R$-module isomorphism
$$
\xymatrix{
f_{iso}:R_{\bf y} \ar[rrr]^{y_1^{n_1}y_2^{n_2}\cdots y_m^{n_m} + p} &&& R_{\geq }.
}
$$
Moreover, there are the natural $R$-module injection
$$
\xymatrix{
f_{inj}:R_{<} \ar[rrr] &&& R_{\bf y},
}
$$
and the natural $R$-module projections
$$
\xymatrix{
f_{proj_{<}}:R_{\bf y} \ar[rrr] &&& R_{<}
},
$$
$$
\xymatrix{
f_{proj_{\geq }}:R_{\bf y} \ar[rrr] &&& R_{\geq }.
}
$$

\indent
The $R$-module $R_{\bf y}^r$ is also split into the direct sum 
$$
R_{\bf y}^r \simeq R_{<}^r \oplus R_{\geq }^r.
$$
Then there are also the $R$-module isomorphism
$$
\xymatrix{
F_{iso}:R_{\bf y}^r \ar[rrr]^{y_1^{n_1}y_2^{n_2}\cdots y_m^{n_m} + p} &&& R_{\geq }^r,
}
$$
the $R$-module injection
$$
\xymatrix{
F_{inj}:R_{<}^r \ar[rrr] &&& R_{\bf y}^r
}
$$
and the $R$-module projections
$$
\xymatrix{
F_{proj_{<}}:R_{\bf y}^r \ar[rrr] &&& R_{<}^r
},
$$

$$
\xymatrix{
F_{proj_{\geq }}:R_{\bf y}^r \ar[rrr] &&& R_{\geq }^r .
}
$$
Since the following $R$-module morphisms $\phi_0$ and $\phi_1$ are $R$-module isomorphisms;
$$
\phi_0 =\xymatrix{
*{\begin{array}{c}
R_{<}^r\\
\oplus\\
R_{\geq }^r\\
\oplus\\
R_{\bf y}^r
\end{array}}
\ar[rrrrrrr]_{\left(
\begin{array}{ccc}
F_{inj}&(y_1^{n_1}y_2^{n_2}\cdots y_m^{n_m} + p)F_{iso}^{-1}&0\\
0&D_0 \, F_{iso}^{-1}&\id_{R_{\bf y}^r}
\end{array}\right)
}&&&&&&&
*{
\begin{array}{c}
R_{\bf y}^r\\
\oplus\\
R_{\bf y}^r ,
\end{array}
}
}$$
$$
\phi_1 =\xymatrix{
*{\begin{array}{c}
R_{<}^r\\
\oplus\\
R_{\geq }^r\\
\oplus\\
R_{\bf y}^r
\end{array}}
\ar[rrrrrrr]_{\left(
\begin{array}{ccc}
F_{inj}&(-y_1^{n_1}y_2^{n_2}\cdots y_m^{n_m} -p)F_{iso}^{-1}&0\\
0&D_1 \, F_{iso}^{-1}&\id_{R_{\bf y}^r}
\end{array}\right)
}&&&&&&&
*{
\begin{array}{c}
R_{\bf y}^r\\
\oplus\\
R_{\bf y}^r ,
\end{array}
}
}$$
the Koszul matrix factorization $K(\mathbf{a};\mathbf{b})_{R_{\bf y}}$ is isomorphic to the following matrix factorization $\overline{M}$ by the isomorphism $\overline{\phi}=$($\phi_0$,$\phi_1$),
\begin{equation*}
\overline{M}=
\left(
\xymatrix{
*{\begin{array}{c}
R_{<}^r\\
\oplus\\
R_{\geq }^r\\
\oplus\\
R_{\bf y}^r
\end{array}}
\ar[rrr]_{\phi_1^{-1}\, D_0 \acute{} \, \phi_0}&&&
*{\begin{array}{c}
R_{<}^r\\
\oplus\\
R_{\geq }^r\\
\oplus\\
R_{\bf y}^r
\end{array}}\ar[rrr]_{\phi_0^{-1}\, D_1 \acute{} \, \phi_1}&&&
*{\begin{array}{c}
R_{<}^r\\
\oplus\\
R_{\geq }^r\\
\oplus\\
R_{\bf y}^r
\end{array}}
}
\right)
.
\end{equation*}
Since $\overline{\phi}^{-1}$ consists of 
$$\phi_0^{-1} = \left(
\begin{array}{cc}
F_{proj_{<}}&0\\
(y_1^{n_1}y_2^{n_2}\cdots y_m^{n_m} + p) F_{iso}^{-1}\, F_{proj_{\geq }}&0\\
-D_0 F_{iso}^{-1}\, F_{proj_{\geq }}&\id_{R_{\bf y}}^r
\end{array}
\right) {\rm and \,\,} \phi_1^{-1} = \left(
\begin{array}{cc}
F_{proj_{<}}&0\\
(-y_1^{n_1}y_2^{n_2}\cdots y_m^{n_m} -p) F_{iso}^{-1}\, F_{proj_{\geq }}&0\\
D_1 F_{iso}^{-1}\, F_{proj_{\geq }}&\id_{R_{\bf y}}^r
\end{array}
\right) ,$$
the morphisms $\phi_1^{-1}\, D_0 \acute{} \, \phi_0$ and $\phi_0^{-1}\, D_1 \acute{} \, \phi_1$ are described by
$$
\phi_1^{-1}\, D_0 \acute{} \, \phi_0=\left(
\begin{array}{ccc}
F_{proj_{<}}\, D_0\, F_{inj} & 0 &0\\
G_1& 0 &F_{iso}\\
G_2& \omega F_{iso}^{-1} &0
\end{array}
\right) {\rm and} \,\,
\phi_0^{-1}\, D_1 \acute{} \, \phi_1=\left(
\begin{array}{ccc}
F_{proj_{<}}\, D_1\, F_{inj} & 0 &0\\
H_1& 0 &F_{iso}\\
H_2& \omega F_{iso}^{-1} &0
\end{array}
\right).
$$
\indent By the construction we have
\begin{eqnarray*}
G_1\,F_{proj_{<}}\, D_1\, F_{inj}+F_{iso}\,H_2&=&0,\\
G_2\,F_{proj_{<}}\, D_1\, F_{inj}+\omega F_{iso}^{-1}\, H_1&=&0,\\
H_1\phi_1^{-1}\, D_0 \acute{} \, \phi_0+F_{iso}G_2&=&0,\\
H_2\phi_1^{-1}\, D_0 \acute{} \, \phi_0+\omega F_{iso}^{-1}G_1&=&0.
\end{eqnarray*}
By the isomorphism $\overline{\psi}=
\left(\left(
\begin{array}{ccc}\id_{R_{<}^r}&0&0\\0&\id_{R_{\geq}^r}&0\\F_{iso}^{-1}G_1&0&\id_{R_{\bf y}^r}\end{array}
\right),\left(
\begin{array}{ccc}\id_{R_{<}^r}&0&0\\0&\id_{R_{\geq}^r}&0\\F_{iso}^{-1}H_1&0&\id_{R_{\bf y}^r}\end{array}
\right)\right)$, the matrix factorization $\overline{M}$ is isomorphic to the following matrix factorization,
\begin{equation*}
\overline{M}=
\left(
\xymatrix{
*{\begin{array}{c}
R_{<}^r\\
\oplus\\
R_{\geq }^r\\
\oplus\\
R_{\bf y}^r
\end{array}}
\ar[rrrrr]_{\left(
\begin{array}{ccc}
\phi_1^{-1}\, D_0 \acute{} \, \phi_0&0&0\\
0&0&F_{iso}\\
0&\omega F_{iso}^{-1}&0
\end{array}
\right)
}&&&&&
*{\begin{array}{c}
R_{<}^r\\
\oplus\\
R_{\geq }^r\\
\oplus\\
R_{\bf y}^r
\end{array}}\ar[rrrrr]_{\left(
\begin{array}{ccc}
\phi_1^{-1}\, D_1 \acute{} \, \phi_0&0&0\\
0&0&F_{iso}\\
0&\omega F_{iso}^{-1}&0
\end{array}
\right)
}&&&&&
*{\begin{array}{c}
R_{<}^r\\
\oplus\\
R_{\geq }^r\\
\oplus\\
R_{\bf y}^r
\end{array}}
}
\right)
.
\end{equation*}
The submatrix factorization
$$
\left(
\xymatrix{
*^{
\left(
\begin{array}{c}
R_{\geq }^r \\
\oplus \\
R_{\bf y}^r
\end{array}
\right)
}
\ar[rrr]_{
\left(
\begin{array}{cc}
0 & F_{iso}\\
\omega F_{iso}^{-1} & 0
\end{array}
\right)
}
&&&
*^{
\left(
\begin{array}{c}
R_{\geq }^r \\
\oplus \\
R_{\bf y}^r
\end{array}
\right)
}
\ar[rrr]_{
\left(
\begin{array}{cc}
0 & F_{iso}\\
\omega F_{iso}^{-1} & 0
\end{array}
\right)
} &&&
*^{
\left(
\begin{array}{c}
R_{\geq }^r \\
\oplus \\
R_{\bf y}^r
\end{array}
\right)
}}
\right)
$$
is contractible in $\HMF^{gr}_{R,\omega }$.
Thus, $\overline{M}$ is isomorphic to the following matrix factorization
$$
\left(
\xymatrix{
*{R_{<}^r}\ar[rrr]^{F_{proj_{<}}\, D_0 \, F_{inj}}&&&*{R_{<}^r}\ar[rrr]^{F_{proj_{<}}\, D_1 \, F_{inj}}&&&*{R_{<}^r}.
}
\right)
$$ in $\HMF^{gr}_{R,\omega }$. 
By the choice of a basis of $R_{\bf y}$ as an $R$-module, we find that this matrix factorization is isomorphic to
 $K(\stackrel{j}{\check{\mathbf{a}}};\stackrel{j}{\check{\mathbf{b}}})_{R_{\bf y}/\left< b_j \right> }$.
\end{proof}

\begin{cor}\label{cor2-10}%corollary
We put $R=\Q [\underline{x}]$ and $R_{\bf y} = R [\underline{y}]$, where $\underline{x}=(x_1,x_2,\cdots ,x_l)$ and $\underline{y}=(y_1,y_2,\cdots ,y_m)$. 
If $\mathbf{a}={}^t( a_1 , a_2 , \cdots , a_k )$ and $\mathbf{b}={}^t( b_1 , b_2 , \cdots , b_k )$, where $a_i ,b_i \in R_{\bf y}$, satisfy the following conditions,
\begin{itemize}
\item[(i)] $\sum_{i=1}^{k}a_i b_i \,(=:\omega) \in R$,\\
\item[($\ast$)] There exist a polynomial $b_j\in R_{\bf y}$ such that ,when $b_j(\underline{x},\underline{y})$ denotes $b_j(\underline{0},\underline{y}) + p$ ($p\in R_{\bf y}$), every monomial of $p$ is not zero in the quotient ring $R_{\bf y}/\left< b_j(\underline{0},\underline{y}) \right>$,
\end{itemize}
then there is the following isomorphism in $\HMF^{gr}_{R,\omega }$,
$$
K(\mathbf{a};\mathbf{b})_{R_{\bf y}} \simeq K(\stackrel{j}{\check{\mathbf{a}}};\stackrel{j}{\check{\mathbf{b}}})_{R_{\bf y}/\left< b_j\right> }.
$$
\end{cor}
\begin{proof}
Each monomial of $b_j(\underline{0},\underline{y})$ is suitable for $y_1^{i_1}y_2^{i_2}\cdots y_m^{i_m}$ of the condition (ii) of Theorem \ref{exclude}. Thus, we obtain this corollary.
\end{proof}

\begin{cor}\label{cor2-11}%corollary
We put $R=\Q [\underline{x}]$ and $R_{\bf y} = R [\underline{y}]$, where $\underline{x}=(x_1,x_2,\cdots ,x_l)$ and $\underline{y}=(y_1,y_2,\cdots ,y_m)$. 
If $\mathbf{a}={}^t( a_1 , a_2 , \cdots , a_k )$ and $\mathbf{b}={}^t( b_1 , b_2 , \cdots , b_k )$, where $a_i ,b_i \in R_{\bf y}$, satisfy the following conditions,
\begin{itemize}
\item[(i)] $\sum_{i=1}^{k}a_i b_i \,(=:\omega) \in R$,\\
\item[(ii)] There are polynomials $b_{j_1},b_{j_2},\cdots ,b_{j_l}$ holding the condition {\rm ($\ast$)} of Corollary \ref{cor2-10} and the condition that the sequence $$(b_{j_1}(\underline{0},\underline{y}),b_{j_2}(\underline{0},\underline{y}),\cdots ,b_{j_l}(\underline{0},\underline{y}))$$ forms regular,
\end{itemize}
then there is the following isomorphism in $\HMF^{gr}_{R,\omega }$,
$$
K(\mathbf{a};\mathbf{b})_{R_{\bf y}} 
\simeq K(\stackrel{j_1,j_2,\cdots ,j_l}{\check{\mathbf{a}}};\stackrel{j_1,j_2,\cdots ,j_l}{\check{\mathbf{b}}})_{R_{\bf y}/\left< b_{j_1},b_{j_2},\cdots ,b_{j_l}\right> }.
$$
\end{cor}

\begin{proof}
Since the sequence $(b_{j_1}(\underline{0},\underline{y}),b_{j_2}(\underline{0},\underline{y}),\cdots ,b_{j_l}(\underline{0},\underline{y}))$ is regular. Then the sequences $\stackrel{j_l}{\check{\mathbf{a}}}$ and $\stackrel{j_l}{\check{\mathbf{b}}}$ satisfy the conditions (i) and (ii) after we apply Corollary \ref{cor2-10} to the polynomial $b_{j_l}(\underline{0},\underline{y})$. Thus we can prove this corollary using induction.
\end{proof}

%%%%%%%%%%%%%%%%%%%%%%%%%%%%%%%%%%%%%%%%%%%%%%%%%%%%%%%%%%%%%%%%%%%%%%%%%%%%%%%%%%%%%%%%%%%%%%%%%%%%%%%%
%
%
% Section3 Matrix factorizations for colored planar diagrams
%
%
%%%%%%%%%%%%%%%%%%%%%%%%%%%%%%%%%%%%%%%%%%%%%%%%%%%%%%%%%%%%%%%%%%%%%%%%%%%%%%%%%%%%%%%%%%%%%%%%%%%%%%%%
\section{Matrix factorizations for colored planar diagrams}\label{sec4}
%%%%%%%%%%%%%%%%%%%%%%%%%%%%%%%%%%%%%%%%%%%%%%%%%%%%%%%%%%%%%%%%%%%%%%%%%%%%%%%%%%%%%%%%%%%%%%%%%%%%%%%%
%
%
% subsection3-1{Poincar\'e series of a graded algebra}
%
%
%%%%%%%%%%%%%%%%%%%%%%%%%%%%%%%%%%%%%%%%%%%%%%%%%%%%%%%%%%%%%%%%%%%%%%%%%%%%%%%%%%%%%%%%%%%%%%%%%%%%%%%%
\subsection{Poincar\'e series of a $\Z$-graded algebra}\label{Poincare}
We briefly recall the Poincar\'e series of a $\Z$-graded algebra. %(For instance, see Section 2.3. in \cite{BT}). 
We only state its definition and properties.\\
\indent
Let $\displaystyle R=\bigoplus_{i=0}^{\infty} R_i$ be a $\Z$-graded algebra over $\Q$. The {\bf Poincar\'e series} of $R$ is defined to be
$$
P_q(R):=\sum_{i=0}^{\infty} (\dim_{\Q} R_i) q^i.
$$
If $x$ is a homogeneous element with $n$-grading in $R$ and not a zero divisor, then we find
$$
P_q(R/xR)=P_q(R)(1-q^n).
$$
Generally, for a regular sequence $\mathcal{X}_r=$($x_1,x_2,\cdots ,x_r$) whose each variable $x_i$ has homogeneous $n_i$-gradings respectively, we find
\begin{equation}
\label{poincare1} P_q(R/\left< \mathcal{X}_r \right>)=P_q(R)(1-q^{n_1})(1-q^{n_2})\cdots (1-q^{n_r}),
\end{equation}
where $\left< \mathcal{X}_r \right>$ is an ideal generated by $x_1,x_2,\cdots ,x_r$.
\begin{pro}\label{poincare3}%proposition
The homogeneous $\Z$-grading terms of 
$$
(1+ x_1+x_2+\cdots +x_r )(1+ y_1+y_2+\cdots +y_s)-1,
$$
where $\deg \, x_i =2i$ and $\deg \, y_i =2i$, form a regular sequence in $\Q [x_1,x_2,\cdots ,x_r,y_1,y_2,\cdots ,y_s]$.
\end{pro}
Let $I$ be an ideal generated by homogeneous $\Z$-grading terms of 
$$
(1+ x_1+x_2+\cdots +x_r )(1+ y_1+y_2+\cdots +y_s)-1,
$$ 
Since these homogeneous terms form a regular sequence, by this proposition and the equation (\ref{poincare1}) the Poincar\'e series of $R/I$ is
\begin{equation}
\label{poincare2}P_q(R/I)=\frac{(1-q^2)(1-q^4)\cdots (1-q^{2r+2s})}{(1-q^2)(1-q^4)\cdots (1-q^{2r})(1-q^2)(1-q^4)\cdots (1-q^{2s})}.
\end{equation}
\indent
Let $x_i$ be the elementary symmetric functions in $r$ variables $t_1,t_2,\cdots ,t_r$ (a $\Z$-grading of every variable $t_i$ is $2$);
\begin{equation}
x_i = \sum_{1\leq n_1 < n_2 < \cdots < n_i \leq r} t_{n_1}t_{n_2}\cdots t_{n_i}
\end{equation}
and let $F_r$ be an $i$ variables function obtained by expanding the power sum $t_1^{n+1}+t_2^{n+1}+\cdots +t_r^{n+1}$ with 
the elementary symmetric functions $x_1,x_2,\cdots ,x_r$, i.e.
\begin{equation}
\label{pow}F_i(x_{1} ,x_{2} ,\cdots ,x_{i})=t_{1}^{n+1} + t_{2}^{n+1} + \cdots + t_{i}^{n+1}.
\end{equation}
The $\Z$-grading of $x_i$ is $2i$ and one of $F_r$ is $2n+2$.
Using a sequence of $i$ variables $\mathcal{X}_{i}=(x_{1} ,x_{2} ,\cdots ,x_{i})$ We often denotes $F_i(x_{1} ,x_{2} ,\cdots ,x_{i})$ by $F_i(\mathcal{X}_{i})$.\\
\indent
We define {\bf Jacobi algebra} $J_{F_r}$ to be
$$
J_{F_r}=\Q [x_1,x_2,\cdots ,x_r] / \left< \frac{\partial F_r}{\partial x_1}, \frac{\partial F_r}{\partial x_2},\cdots ,\frac{\partial F_r}{\partial x_r} \right>.
$$ 
\indent
In the work of Gepner \cite{Gepner}, the Jacobi algebra $J_{F_r}$ is isomorphic to the cohomology ring of the Grassmannian manifold $\H^* (Gr_r (\C^n))$ as a $\Z$-graded algebra. 
Since the Poincar\'e series of the cohomology ring $\H^* (Gr_r (\C^n))$ is
$$
\frac{(1-q^2)(1-q^4)\cdots (1-q^{2n})}{(1-q^2)(1-q^4)\cdots (1-q^{2r})(1-q^2)(1-q^4)\cdots (1-q^{2n-2r})},
$$
the Poincar\'e series of Jacobi algebra $J_{F_r}$ is also
\begin{equation}
\label{Jacobi} P_q(J_{F_r})=\frac{(1-q^2)(1-q^4)\cdots (1-q^{2n})}{(1-q^2)(1-q^4)\cdots (1-q^{2r})(1-q^2)(1-q^4)\cdots (1-q^{2n-2r})}.
\end{equation}
%%%%%%%%%%%%%%%%%%%%%%%%%%%%%%%%%%%%%%%%%%%%%%%%%%%%%%%%%%%%%%%%%%%%%%%%%%%%%%%%%%%%%%%%%%%%%%%%%%%%%%%%
%
%
% subsection3-2 Colored planar diagrams and graded matrix factorizations
%
%
%%%%%%%%%%%%%%%%%%%%%%%%%%%%%%%%%%%%%%%%%%%%%%%%%%%%%%%%%%%%%%%%%%%%%%%%%%%%%%%%%%%%%%%%%%%%%%%%%%%%%%%%
\subsection{Colored planar diagrams and matrix factorizations}\label{sec4.2}
We define matrix factorizations for the colored planar diagrams in Figure $5$.
Let $t_{i,k}$ be a variable with $\Z$-grading $2$ for any formal indexes $i$ and $k$, let $x_{j,k}$ be the elementary symmetric function with $\Z$-grading $2j$ in a polynomial ring $\Q [t_{1,k},t_{2,k},\cdots ,t_{i,k}]$, where $k$ is a formal index, 
and let $F_i$ be an $i$ variables function obtained by expanding the power sum $t_{1,k}^{n+1} + t_{2,k}^{n+1} + \cdots + t_{i,k}^{n+1}$ 
with the elementary symmetric function $x_{1,k} ,x_{2,k} ,\cdots ,x_{i,k}$ (as Equation (\ref{pow})).

\indent
We assign a sequence $\mathcal{X}_{i,k}$ of $i$ variables $x_{1,k}$, $x_{2,k}$, $\cdots$, $x_{i,k}$ on a boundary of a line colored $i$ and define a map $\c\acute{}$ from a colored planar diagram with such an assignment to a matrix factorization.
\begin{de}
\indent
A matrix factorization for 
\begin{center}
\input{figmoy3-mf}\hspace{2cm}${\rm (i=1,\cdots ,n )}$
\end{center}
\hspace{1mm}\\\\
is defined to be
\begin{equation}
\label{line-mf}\c\acute{} \Biggl( \input{figmoy3-mf} \Biggl)_n:=
\mathop{\boxtimes}_{j=1}^{i} K\Big( L^{1;2}_{j,i} ;x_{j,1}-x_{j,2} \Big)_{\Q [\mathcal{X}_{i,1} ,\mathcal{X}_{i,2} ]} ,
\end{equation}
where 
$$
L^{1;2}_{j,i} =\frac{F_{i}(x_{1,2},\cdots ,x_{j-1,2} ,x_{j,1},x_{j+1,1},\cdots ,x_{i,1})-F_{i}(x_{1,2},\cdots ,x_{j-1,2},x_{j,2},x_{j+1,1},\cdots ,x_{i,1})}{x_{j,1}-x_{j,2}}
$$
and 
$$
\Q [\mathcal{X}_{i,1} ,\mathcal{X}_{i,2} ] = \Q [x_{1,1},\cdots ,x_{i,1},x_{1,2},\cdots ,x_{i,2} ].
$$
\end{de}
This matrix factorization is an object of $\MF^{gr}_{\Q [\mathcal{X}_{i,1} ,\mathcal{X}_{i,2} ],F_i (\mathcal{X}_{i,1}) -F_i (\mathcal{X}_{i,2}) }$.
\begin{rem}
For $i \geq n+1$, we can also consider the matrix factorization for a line colored $i$, \input{figmoy5} as the above definition.
However, we find that these matrix factorizations are homotopic to the zero matrix factorization.
\end{rem}
\begin{de}
\indent
A matrix factorization for the following trivalent diagrams,
\begin{center}
\input{figgluing-in-3valent-mf} and \input{figgluing-out-3valent-mf},\vspace{1cm}
\end{center}
is defined to be
\begin{eqnarray}
\label{n-mf}\c\acute{} \left(\input{figgluing-in-3valent-mf} \right)_n &:=& \mathop{\boxtimes}_{j=1}^{i_{3}} 
K\Big( \Lambda_{j,i_1;i_2}^{1;2,3} ;x_{j,3}-X^{1;2}_{j,i_1;i_2} \Big)_{\Q [\mathcal{X}_{i_{1},1},\mathcal{X}_{i_{2},2},\mathcal{X}_{i_{3},3}]} ,\\[-0.1em]
\label{v-mf}\c\acute{} \left(\input{figgluing-out-3valent-mf} \right)_n &:=& \mathop{\boxtimes}_{j=1}^{i_{3}} 
K\Big( V_{j,i_1;i_2}^{1;2,3} ;X^{1;2}_{j,i_1;i_2}-x_{j,3} \Big)_{\Q [\mathcal{X}_{i_{1},1},\mathcal{X}_{i_{2},2},\mathcal{X}_{i_{3},3}]} \{ - i_1 i_2 \} ,
\end{eqnarray}
where 
\begin{eqnarray*}
\Lambda_{j,i_1;i_2}^{1;2,3}=\frac{F_{i_{3}}(X^{1;2}_{1,i_1;i_2},\cdots ,X^{1;2}_{j-1,i_1;i_2} ,x_{j,3},x_{j+1,3},\cdots ,x_{i_{3},3})-F_{i_{3}}(X^{1;2}_{1,i_1;i_2},\cdots ,X^{1;2}_{j-1,i_1;i_2},X^{1;2}_{j,i_1;i_2},x_{j+1,3},\cdots ,x_{i_{3},3})}{x_{j,3}-X^{1;2}_{j,i_1;i_2}},\\[-0.1em]
V_{j,i_1;i_2}^{1;2,3}=\frac{F_{i_{3}}(x_{1,3},\cdots ,x_{j-1,3} ,X^{1;2}_{j,i_1;i_2},X^{1;2}_{j+1,i_1;i_2},\cdots ,X^{1;2}_{i_3,i_1;i_2})-F_{i_{3}}(x_{1,3},\cdots ,x_{j-1,3},x_{j,3},X^{1;2}_{j+1,i_1;i_2},\cdots ,X^{1;2}_{i_3,i_1;i_2})}{X^{1;2}_{j,i_1;i_2}-x_{j,3}},
\end{eqnarray*}
and $X^{k_1;k_2}_{j,i_1;i_2}$ is the $2j$-grading term of
$(1+ x_{1,k_1}+x_{2,k_1}+\cdots +x_{i_1,k_1} )(1+ x_{1,k_2}+x_{2,k_2}+\cdots +x_{i_2,k_2} )-1$. \\
\end{de}
Two matrix factorizations (\ref{n-mf}) and (\ref{v-mf}) are an object of $\MF^{gr}_{\Q [\mathcal{X}_{i_1,1} ,\mathcal{X}_{i_2,2},\mathcal{X}_{i_3,3} ],F_{i_3} (\mathcal{X}_{i_3,3}) -F_{i_1} (\mathcal{X}_{i_1,1}) -F_{i_2}(\mathcal{X}_{i_2,2}) }$ and an object of $\MF^{gr}_{\Q [\mathcal{X}_{i_1,1} ,\mathcal{X}_{i_2,2},\mathcal{X}_{i_3,3} ], F_{i_1} (\mathcal{X}_{i_1,1}) +F_{i_2} (\mathcal{X}_{i_2,2}) - F_{i_3} (\mathcal{X}_{i_3,3})}$ respectively.\\
\indent
We define a matrix factorization for a more general planar diagram by using tensor product as follows. \\
\indent
We consider two planar diagrams which have a line colored $i$ 
and can be match with keeping the orientation on the line colored $i$, 
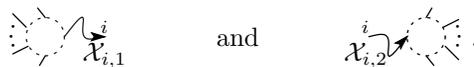
\begin{center}
\input{figgluing-gen-mf1} \hspace{1cm}and\hspace{1cm} \input{figgluing-gen-mf2} .
\end{center}
Using the matrix factorizations 
$$
\c\acute{}\Big(\input{figgluing-gen-mf1}\Big)_{n} \in \ob (\MF^{gr}_{R,\omega + F_{i}(\mathcal{X}_{i,1})})
\hspace{0.5cm} {\rm and} \hspace{0.5cm} 
\c\acute{}\Big(\input{figgluing-gen-mf2}\Big)_{n}
\in \ob (\MF^{gr}_{R\acute{},\omega\acute{} - F_{i}(\mathcal{X}_{i,2})}),
$$
we define a matrix factorization for the glued planar diagram
$$
\c\acute{}\Bigg( \input{figgluing-gen2} \Bigg)_{n} \in \ob (\MF^{gr}_{R \ostimes R\acute{},\omega + \omega\acute{}}) 
\hspace{0.5cm} {\rm by} \hspace{0.5cm} 
\c\acute{}\Big(\input{figgluing-gen-mf1}\Big)_{n} \boxtimes \c\acute{}\Big(\input{figgluing-gen-mf2}\Big)_{n}\Big|_{\mathcal{X}_{i,2}\to\mathcal{X}_{i,1}}.
$$
This means that we identify the sequence $\mathcal{X}_{i,1}$ and the sequence $\mathcal{X}_{i,2}$ after taking the tensor product of these matrix factorizations. Since a potential of the tensor product of two matrix factorizations is the sum of each potential, the potential of the glued matrix factorization is $\omega+\omega\acute{}$.\\
\indent For a planar diagram $\Gamma$ composed of the disjoint union of planar diagrams $\Gamma_1$ and $\Gamma_2$, we define
$$
\c\acute{}(\ \Gamma\ )_{n} := \c\acute{}(\ \Gamma_1\ )_{n}\boxtimes \c\acute{}(\ \Gamma_2\ )_{n}
$$

Since we only consider a colored planar diagram decomposing into the colored trivalent diagrams in Figure $5$, we can obtain a matrix factorization for the colored planar diagram by taking tensor product of matrix factorizations for the trivalent diagrams.
By definition of gluing, the potential of a colored planar diagram depend only on the boundary assignment of colored planar diagram.
\begin{pro}%proposition
{\rm\bf (1)}There is the following isomorphism in $\HMF^{gr}_{\Q [\mathcal{X}_{i,1}]\ostimes R,\omega +F_i(\mathcal{X}_{i,1})}$, where the polynomial ring $R$ 
and the potential $\omega$ are determined by the boundary sequences except the sequence $\mathcal{X}_{i,1}$,

\begin{eqnarray*}
\c\acute{}\Big( \ \input{fig-gen-gluing-mf3} \ \Big)_{n} \simeq \c\acute{}\Big( \ \input{fig-gen-gluing-mf1} \ \Big)_{n}
\end{eqnarray*}

{\rm\bf (2)}There is the following isomorphism in $\HMF^{gr}_{\Q [\mathcal{X}_{i,1}]\ostimes R,\omega -F_i(\mathcal{X}_{i,1})}$, where the polynomial ring $R$ 
and the potential $\omega$ are determined by the boundary sequences except the sequence $\mathcal{X}_{i,1}$,

\begin{eqnarray*}
\c\acute{}\Big( \ \input{fig-gen-gluing-mf4} \ \Big)_{n} \simeq \c\acute{}\Big( \ \input{fig-gen-gluing-mf2} \ \Big)_{n}
\end{eqnarray*}
\\

{\rm\bf (3)}There is the following isomorphism in $\HMF^{gr}_{\Q ,0}$

\begin{eqnarray*}
\c\acute{}\Big( \ \input{figcircle-weight-r-mf} \ \Big)_{n} &=& \mathop{\boxtimes}_{j=1}^{i} K\Big( L^{1;2}_{j,i} ;x_{j,1}-x_{j,2} \Big)_{\Q [\mathcal{X}_{i,1} ,\mathcal{X}_{i,2} ]}\Big|_{\mathcal{X}_{i,2} \to \mathcal{X}_{i,1}}\\[-0.1em]
&\simeq& (\, J_{F_{i}(\mathcal{X}_{i,1})} \to 0\to J_{F_{i}(\mathcal{X}_{i,1})}  \,)\left\{ -in+i^2 \right\} \left< i \right>,
\end{eqnarray*}

where $J_{F_{i}(\mathcal{X}_{i,1})}$ is Jacobi algebra for the polynomial $F_{i}(\mathcal{X}_{i,1}))$, i.e.
$$
J_{F_{i}(\mathcal{X}_{i,1})}=\Q [\mathcal{X}_{i,1}] \left/ 
\left< \frac{\partial F_{i}}{\partial x_{1,1}}, \cdots , \frac{\partial F_{i}}{\partial x_{i,1}} \right>\right. .
$$
\end{pro}

\begin{proof}

{\bf (1)}We describe the matrix factorization $\c\acute{}\Big( \ \input{fig-gen-gluing-mf5} \ \Big)_{n}$ in $ \HMF^{gr}_{\Q [\mathcal{X}_{i,2}]\ostimes R,\omega +F_i(\mathcal{X}_{i,2})}$ as $\left( M_0,M_1,D_0,D_1 \right) $,
where $M_j$ is an $R \ostimes \Q [\mathcal{X}_{i,2}]$-module. Then we have
$$
\c\acute{}\Big( \ \input{fig-gen-gluing-mf3} \ \Big)_{n} = \left( M_0,M_1,D_0,D_1 \right)
\mathop{\boxtimes}_{j=1}^{i} K\left( L_{j,i}^{1;2};x_{j,1}-x_{j,2} \right)_{\Q [\mathcal{X}_{i,1} ,\mathcal{X}_{i,2} ]}.
$$
We can regard this matrix factorization composed of $\Q [\mathcal{X}_{i,1},\mathcal{X}_{i,2}] \ostimes R$-modules as 
a matrix factorization which consists of infinite rank $\Q [\mathcal{X}_{i,1}] \ostimes R$-modules.
Because the potential of this matrix factorization is in the polynomial ring $\Q [\mathcal{X}_{i,1}] \ostimes R$. 
Moreover, the polynomials  $x_{j,1}-x_{j,2}|_{\mathcal{X}_{i,1}=0}$ ( $j=1,2,\cdots ,i$ ) form regular in $\Q [\mathcal{X}_{i,1},\mathcal{X}_{i,2}] \ostimes R$.
Thus, we can apply Corollary \ref{cor2-11} to $\left( M_0,M_1,D_0,D_1 \right)
\mathop{\boxtimes}_{j=1}^{i} K\left( L_{j,i}^{1;2};x_{j,1}-x_{j,2} \right)_{\Q [\mathcal{X}_{i,1} ,\mathcal{X}_{i,2} ]}$. 
Then we have 
$$
\left( M_0 ,M_1,D_0,D_1\right)\left|_{\mathcal{X}_{i,1}\to\mathcal{X}_{i,2}}\right.
\boxtimes \left(\Q [\mathcal{X}_{i,1},\mathcal{X}_{i,2}] \ostimes R/ \left<x_{1,1}-x_{1,2}, \cdots ,x_{i,1}-x_{i,2} \right>,0,0,0\right)
$$
The quotient polynomial ring 
$\Q [\mathcal{X}_{i,1},\mathcal{X}_{i,2}] \ostimes R/ \left<x_{1,1}-x_{1,2}, \cdots ,x_{i,1}-x_{i,2} \right>$ is isomorphic to
$\Q [\mathcal{X}_{i,1}] \ostimes R$.
Thus, we obtain the right-hand matrix factorization.\\

{\bf (2)}This proof is similar to the proof of (1).\\

{\bf (3)}We have
\begin{equation*}
\mathop{\boxtimes}_{j=1}^{i} 
K\Big( L^{1;2}_{j,i} ;x_{j,1}-x_{j,2} \Big)_{\Q [\mathcal{X}_{i,1} ,\mathcal{X}_{i,2} ]}\Big|_{\mathcal{X}_{i,2} \to \mathcal{X}_{i,1}} =
\mathop{\boxtimes}_{j=1}^{i} 
\left(
	\Q [\mathcal{X}_{i,1}],\Q [\mathcal{X}_{i,1}]\{ 2j-1-n \},L^{1;2}_{j,i}|_{\mathcal{X}_{i,2} \to \mathcal{X}_{i,1}}0 \right).
\end{equation*}
The polynomial $L^{1;2}_{j,i}|_{\mathcal{X}_{i,2} \to \mathcal{X}_{i,1}}$ is
$$
\frac{F_{i}(\cdots ,x_{j-1,2},x_{j,1},x_{j+1,1},\cdots )-F_{i}(\cdots ,x_{j-1,2},x_{j,2},x_{j+1,1},\cdots )}{x_{j,1}-x_{j,2}}
\Big|_{\mathcal{X}_{i,2} \to \mathcal{X}_{i,1}} =\frac{\partial F_{i}(\mathcal{X}_{i,1})}{\partial x_{j,1}}.
$$
Hence, we apply Theorem \ref{exclude} to these polynomials of the matrix factorization after using Proposition \ref{functor1} and \ref{functor2};

\begin{eqnarray*}
&&\mathop{\boxtimes}_{j=1}^{i} 
K\Big( L^{1;2}_{j,i} ;x_{j,1}-x_{j,2} \Big)_{\Q [\mathcal{X}_{i,1} ,\mathcal{X}_{i,2} ]}\Big|_{\mathcal{X}_{i,2}\to\mathcal{X}_{i,1}}\\
&&\\[-0.1em]
&\simeq&
\mathop{\boxtimes}_{j=1}^{i} 
\left( \Q [\mathcal{X}_{i,1}],\Q [\mathcal{X}_{i,1}]\{ n+1-2j \},0,\frac{\partial F_{i}(\mathcal{X}_{i,1})}{\partial x_{j,1}}\right) \{ -i n + i^2 \} \left< i \right> \\[-0.1em]
&&\\[-0.1em]
&\simeq&\left(J_{F_{i}(\mathcal{X}_{i,1})},0,0,0\right)
\left\{ -in+i^2 \right\} \left< i \right>.
\end{eqnarray*}

\end{proof}

A matrix factorization for the loop colored $i$, %WinTpicVersion2.15
\unitlength 0.1in
\begin{picture}(3.0,1.50)(2.2500,-4.0)
% CIRCLE 2 0 3 0
% 4 350 750 350 900 350 900 350 900
% 
\special{pn 8}%
\special{ar 350 350 120 120  0.0000000 6.2831853}%
% STR 2 0 3 0
% 3 1770 850 1770 950 5 0
% ${}_1$
\put(5.25000,-3.5000){\makebox(0,0){${}_{i}$}}%
\end{picture}%
, is defined to be the above matrix factorization;
$$
\MF\acute{}\left( %WinTpicVersion2.15
\unitlength 0.1in
\begin{picture}(3.0,1.50)(2.2500,-4.0)
% CIRCLE 2 0 3 0
% 4 350 750 350 900 350 900 350 900
% 
\special{pn 8}%
\special{ar 350 350 120 120  0.0000000 6.2831853}%
% STR 2 0 3 0
% 3 1770 850 1770 950 5 0
% ${}_1$
\put(5.25000,-3.5000){\makebox(0,0){${}_{i}$}}%
\end{picture}%
 \right):=\left(J_{F_{i}(\mathcal{X}_{i,1})},0,0,0\right)
\left\{ -in+i^2 \right\} \left< i \right>.
$$
\begin{cor}%corollary
The Euler characteristic of the homology $\H (\c\acute{} (%WinTpicVersion2.15
\unitlength 0.1in
\begin{picture}(3.0,1.50)(2.2500,-4.0)
% CIRCLE 2 0 3 0
% 4 350 750 350 900 350 900 350 900
% 
\special{pn 8}%
\special{ar 350 350 120 120  0.0000000 6.2831853}%
% STR 2 0 3 0
% 3 1770 850 1770 950 5 0
% ${}_1$
\put(5.25000,-3.5000){\makebox(0,0){${}_{i}$}}%
\end{picture}%
)_n)$ equals MOY polynomial
for the loop colored $i$;
$$
\overline{\chi}\Big( \H (\c\acute{}\Big(%WinTpicVersion2.15
\unitlength 0.1in
\begin{picture}(3.0,1.50)(2.2500,-4.0)
% CIRCLE 2 0 3 0
% 4 350 750 350 900 350 900 350 900
% 
\special{pn 8}%
\special{ar 350 350 120 120  0.0000000 6.2831853}%
% STR 2 0 3 0
% 3 1770 850 1770 950 5 0
% ${}_1$
\put(5.25000,-3.5000){\makebox(0,0){${}_{i}$}}%
\end{picture}%
\Big)_n) \Big) =\left[ n \atop i \right].
$$
\end{cor}
\begin{proof}
By definition of Euler characteristic, the left-hand side of this equation equals the Poincar\'e series of Jacobi algebra times $q^{-in+i^2}$.
Thus, we obtain the right-hand side of this equation by using the equation (\ref{Jacobi}) in Section \ref{Poincare}.
\end{proof}
\begin{pro}\label{mat-equiv1}%proposition
Let $R$ be a polynomial ring generated by variables of sequences $\mathcal{X}_{i,1}$, $\mathcal{X}_{i,2}$, $\mathcal{X}_{i,3}$, $\mathcal{X}_{i,4}$ 
and let $\omega$ be a polynomial
$$
F_{i_4}(x_{1,4},\cdots ,x_{i_4,4})-F_{i_1}(x_{1,1},\cdots ,x_{i_1,1})-F_{i_2}(x_{1,2},\cdots ,x_{i_2,2})-F_{i_3}(x_{1,3},\cdots ,x_{i_3,3}).
$$

\begin{enumerate}
\item[{\rm\bf (1)}] There is the following isomorphism in $\HMF^{gr}_{R,\omega}$
$$
\c\acute{}\left( \input{fig-ass-dia1-mf}\right)_n \simeq \c\acute{}\left( \input{fig-ass-dia2-mf}\right)_n ,
$$ 
\item[{\rm\bf (2)}] There is the following isomorphism in $\HMF^{gr}_{R,-\omega}$
$$
\c\acute{}\left( \input{fig-coass-dia1-mf}\right)_n \simeq \c\acute{}\left( \input{fig-coass-dia2-mf}\right)_n ,
$$
\end{enumerate}

where $1 \leq i_1,i_2,i_3 \leq n-2$, $i_5 = i_1 +i_2\leq n-1$, $i_6 = i_2 + i_3\leq n-1$ and $i_4 = i_1 + i_2 + i_3 \leq n$.
\end{pro}

\begin{proof}
{\bf (1)}
We have
\begin{eqnarray*}
&&\c\acute{}\left( \input{fig-ass-dia1-mf}\right)_n \\[-0.1em]
&=& \mathop{\boxtimes}_{j=1}^{i_4} 
K\Big( \Lambda_{j,i_5;i_3}^{5;3,4} ;x_{j,4}-X^{5;3}_{j,i_5;i_3} \Big)_{\Q [\mathcal{X}_{i_{3},3},\mathcal{X}_{i_{4},4},\mathcal{X}_{i_{5},5}]} 
\mathop{\boxtimes}_{j=1}^{i_5}
K\Big( \Lambda_{j,i_1;i_2}^{1;2,5} ;x_{j,5}-X^{1;2}_{j,i_1;i_2} \Big)_{\Q [\mathcal{X}_{i_{1},1},\mathcal{X}_{i_{2},2},\mathcal{X}_{i_{5},5}]}.
\end{eqnarray*}
Since the potential of this matrix factorization does not include the variables of $\mathcal{X}_{i_5,5}$ and 
\begin{equation*}
(x_{1,5}-X^{1;2}_{1,i_1;i_2},\cdots ,x_{i_5,5}-X^{1;2}_{i_5,i_1;i_2})|_{(\mathcal{X}_{i_{1},1},\mathcal{X}_{i_{2},2},\mathcal{X}_{i_{3},3},\mathcal{X}_{i_{4},4})=(\underline{0})}=(\mathcal{X}_{i_{5},5})
\end{equation*}
is obviously a regular sequence, we can apply Corollary \ref{cor2-11} to these variables. Then the matrix factorization is isomorphic to
\begin{equation*}
\mathop{\boxtimes}_{j=1}^{i_4} K\Big( \Lambda_{j,i_5;i_3}^{5;3,4} ;x_{j,4}-X^{5;3}_{j,i_5;i_3} \Big)_{\Q [\mathcal{X}_{i_{1},1},\mathcal{X}_{i_{2},2},\mathcal{X}_{i_{3},3},\mathcal{X}_{i_{4},4},\mathcal{X}_{i_{5},5}]/
\left< x_{1,5}-X^{1;2}_{1,i_1;i_2},\cdots ,x_{i_5,5}-X^{1;2}_{i_5,i_1;i_2} \right> }. 
\end{equation*}

In the quotient ring $\Q [\mathcal{X}_{i_{1},1},\mathcal{X}_{i_{2},2},\mathcal{X}_{i_{3},3},\mathcal{X}_{i_{4},4},\mathcal{X}_{i_{5},5}]/\left< x_{1,5}-X^{1;2}_{1,i_1;i_2},\cdots ,x_{i_5,5}-X^{1;2}_{i_5,i_1;i_2}  \right>$, the polynomial $X^{5;3}_{j,i_5;i_3}$ equals the $2j$-grading term of
$$
(1+x_{1,1}+x_{2,1}+\cdots +x_{i_1,1})(1+x_{1,2}+x_{2,2}+\cdots +x_{i_2,2})(1+x_{1,3}+x_{2,3}+\cdots +x_{i_3,3}).
$$
If $X_{j,i_1;i_2;i_3}^{1;2;3}$ denotes the $2j$-grading term of 
$$
(1+x_{1,1}+x_{2,1}+\cdots +x_{i_1,1})(1+x_{1,2}+x_{2,2}+\cdots +x_{i_2,2})(1+x_{1,3}+x_{2,3}+\cdots +x_{i_3,3}),
$$
the polynomial $\Lambda_{j,i_5;i_3}^{5;3,4}$ equals to
$$
\frac{F_{i_4}(\cdots ,X_{j-1,i_1;i_2;i_3}^{1;2;3},x_{j,4},x_{j+1,4},\cdots )-F_{i_4}(\cdots ,X_{j-1,i_1;i_2;i_3}^{1;2;3},X_{j,i_1;i_2;i_3}^{1;2;3},x_{j+1,4},\cdots )}{x_{j,4}-X_{j,i_1;i_2;i_3}^{1;2;3}}.
$$
Let $\Lambda_{j,i_1;i_2;i_3}^{1;2;3,4}$ denote this polynomial.  
Then the matrix factorization is isomorphic to
$$
\mathop{\boxtimes}_{j=1}^{i_4} K\Big( \Lambda_{j,i_1;i_2;i_3}^{1;2;3,4} ;x_{j,4}-X^{1;2;3}_{j,i_1;i_2;i_3} \Big)_{\Q [\mathcal{X}_{i_{1},1},\mathcal{X}_{i_{2},2},\mathcal{X}_{i_{3},3},\mathcal{X}_{i_{4},4}] }. 
$$
\indent
The other one is
\begin{eqnarray*}
&&\c\acute{}\left( \input{fig-ass-dia2-mf}\right)_n \\[-0.1em]
&=& \mathop{\boxtimes}_{j=1}^{i_4} 
K\Big( \Lambda_{j,i_1;i_6}^{1;6,4} ;x_{j,4}-X^{1;6}_{j,i_1;i_6} \Big)_{\Q [\mathcal{X}_{i_{1},1},\mathcal{X}_{i_{4},4},\mathcal{X}_{i_{6},6}]} 
\mathop{\boxtimes}_{j=1}^{i_6}
K\Big( \Lambda_{j,i_2;i_3}^{2;3,6} ;x_{j,6}-X^{2;3}_{j,i_2;i_3} \Big)_{\Q [\mathcal{X}_{i_{2},2},\mathcal{X}_{i_{3},3},\mathcal{X}_{i_{6},6}]} 
\end{eqnarray*}

Since the potential of this matrix factorization does not include the variables of $\mathcal{X}_{i_6,6}$ and
\begin{equation*}
(x_{1,6}-X^{2;3}_{1,i_2;i_3},\cdots ,x_{i_6,6}-X^{2;3}_{i_6,i_2;i_4})|_{(\mathcal{X}_{i_{1},1},\mathcal{X}_{i_{2},2},\mathcal{X}_{i_{3},3},\mathcal{X}_{i_{4},4})=(\underline{0})}=(\mathcal{X}_{i_{6},6})
\end{equation*}
is a regular sequence in $\Q[\mathcal{X}_{i_{1},1}\mathcal{X}_{i_{2},2}\mathcal{X}_{i_{3},3}\mathcal{X}_{i_{4},4}\mathcal{X}_{i_{6},6}]$, we can apply Corollary \ref{cor2-11} to these variables. We similarly obtain the result that the matrix factorization is isomorphic to
$$
\mathop{\boxtimes}_{j=1}^{i_4} K\Big( \Lambda_{j,i_1;i_2;i_3}^{1;2;3,4} ;x_{j,4}-X^{1;2;3}_{j,i_1;i_2;i_3} \Big)_{\Q [\mathcal{X}_{i_{1},1},\mathcal{X}_{i_{2},2},\mathcal{X}_{i_{3},3},\mathcal{X}_{i_{4},4}] }. 
$$
{\bf (2)}The proof of this proposition is similar to (1). By using Corollary \ref{cor2-11}, we find that the left-hand and right-hand matrix factorization are isomorphic to
$$
\mathop{\boxtimes}_{j=1}^{i_4} K\Big( V_{j,i_1;i_2;i_3}^{1;2;3,4} ;X^{1;2;3}_{j,i_1;i_2;i_3}-x_{j,4} \Big)_{\Q [\mathcal{X}_{i_{1},1},\mathcal{X}_{i_{2},2},\mathcal{X}_{i_{3},3},\mathcal{X}_{i_{4},4}] } \left\{ -i_1 i_2-i_2 i_3-i_1 i_3 \right\},
$$
where the polynomial $V_{j,i_1;i_2;i_3}^{1;2;3,4}$ is
$$
\frac{F_{i_4}(\cdots ,x_{j-1,4} ,X_{j,i_1;i_2;i_3}^{1;2;3},X_{j+1,i_1;i_2;i_3}^{1;2;3},\cdots )-F_{i_4}(\cdots ,x_{j-1,4} ,x_{j,4},X_{j+1,i_1;i_2;i_3}^{1;2;3},\cdots )}{X_{j,i_1;i_2;i_3}^{1;2;3}-x_{j,4}}.
$$
\end{proof}

For $0 \leq n_2 \leq j \leq n_1$, we define the positive integer $p(j,n_1;n_2)$ by
$$
\left[ {n_1}\atop{n_2} \right] = \sum_{j=0}^{n_1 n_2-n_2^2} p(j,n_1;n_2) \ q^{-n_1 n_2+n_2^2+2j}.
$$

\begin{pro}\label{mat-equiv2}%proposition
{\rm\bf (1)}There is the following isomorphism in 
$\HMF^{gr}_{\Q [\mathcal{X}_{i_3,1},\mathcal{X}_{i_3,2}],F_{i_3}(\mathcal{X}_{i_3,1})-F_{i_3}(\mathcal{X}_{i_3,2})}$
$$
\c\acute{}\left( \input{fig-bubble-color-mf}\right)_n \simeq 
\bigoplus_{j=0}^{i_1 i_3 -i_1^2}\left( \c\acute{}\left( \input{fig-line-color-mf}\right)_n \{ -i_1 i_3 +i_1^2 +2j \}\right)^{\oplus p(j,i_3,i_1)},
$$
{\rm\bf (2)}There is the following isomorphism in
$\HMF^{gr}_{\Q [\mathcal{X}_{i_1,1},\mathcal{X}_{i_1,2}],F_{i_1}(\mathcal{X}_{i_1,1})-F_{i_1}(\mathcal{X}_{i_1,2})}$
$$
\c\acute{}\left( \input{fig-bubble-color1-mf}\right)_n 
\simeq \bigoplus_{j=0}^{i_2(n-i_1)-i_2^2}\left( \c\acute{}\left( \input{fig-line-color1-mf}\right)_n \{ -i_2(n-i_1)+i_2^2 +2j \}\right)^{\oplus p(j,n-i_1,i_2)}\hspace{-1.5cm}\left< i_2 \right> ,
$$ 
where $1 \leq i_1,i_2 \leq n-1$ and $i_3 = i_1 + i_2 \leq n$.
\end{pro}

\begin{proof}
{\bf(1)}We have
\begin{eqnarray*}
&&\c\acute{}\left( \input{fig-bubble-color-mf}\right)_n \\[-0.1em]
&=& 
\mathop{\boxtimes}^{i_3}_{j=1} K\Big( \Lambda_{j,i_1;i_2}^{3;4,1} ;x_{j,1}+X^{3;4}_{j,i_1;i_2} \Big)_{\Q [\mathcal{X}_{i_{3},1},\mathcal{X}_{i_{1},3},\mathcal{X}_{i_{2},4}] } 
\mathop{\boxtimes}^{i_3}_{j=1} K\Big( V_{j,i_1;i_2}^{3;4,2} ;x_{j,2}+X^{3;4}_{j,i_1;i_2} \Big)_{\Q [\mathcal{X}_{i_{3},2},\mathcal{X}_{i_{1},3},\mathcal{X}_{i_{2},4}] } 
\left\{ -i_1 i_2 \right\} .
\end{eqnarray*}
The potential of this matrix factorization does not include the variables of $\mathcal{X}_{i_1,3}$, $\mathcal{X}_{i_2,4}$ and 
the sequence 
\begin{equation*}
(x_{1,2}+X_{1,i_1;i_2}^{3;4},x_{2,2}+X_{2,i_1;i_2}^{3;4},\cdots ,x_{i_3,2}+X_{i_3,i_1;i_2}^{3;4})
\left|_{(\mathcal{X}_{i_{3},1},\mathcal{X}_{i_{3},2})=(\underline{0})}\right.=(X_{1,i_1;i_2}^{3;4},X_{2,i_1;i_2}^{3;4},\cdots ,X_{i_3,i_1;i_2}^{3;4})
\end{equation*} 
is regular by Proposition \ref{poincare3},
we can apply Corollary \ref{cor2-11} to the matrix factorization. 
Thus we have
\begin{equation*}
\mathop{\boxtimes}^{i_3}_{j=1} K\Big( \Lambda_{j,i_1;i_2}^{3;4,1} ;x_{j,1}-x_{j,2} \Big)_{\Q [\mathcal{X}_{i_{3},1},\mathcal{X}_{i_{3},2},\mathcal{X}_{i_{1},3},\mathcal{X}_{i_{2},4}]/\left< x_{1,2}+X^{3;4}_{1,i_1;i_2},\cdots ,x_{i_3,2}+X^{3;4}_{i_3,i_1;i_2} \right> } \left\{ -i_1 i_2 \right\}.
\end{equation*}

The polynomial $\Lambda_{j,i_1;i_2}^{3;4,1}$ is equal to $L_{j,i_3}^{12}$ in the quotient ring 
$$
\Q [\mathcal{X}_{i_{3},1},\mathcal{X}_{i_{3},2},\mathcal{X}_{i_{1},3},\mathcal{X}_{i_{2},4}]/\left< x_{1,2}+X^{3;4}_{1,i_1;i_2},\cdots ,x_{i_3,2}+X^{3;4}_{i_3,i_1;i_2} \right> .
$$
By the equation (\ref{poincare2}) in Section \ref{Poincare}, the Poincar\'e series of the quotient ring is

\begin{eqnarray*}
&&P_q (\Q [\mathcal{X}_{i_{3},1},\mathcal{X}_{i_{3},2},\mathcal{X}_{i_{1},3},\mathcal{X}_{i_{2},4}]/\left< x_{1,2}+X^{3;4}_{1,i_1;i_2},\cdots ,x_{i_3,2}+X^{3;4}_{i_3,i_1;i_2} \right> ) \\[-0.1em]
&=&\frac{(1-q^2)(1-q^4)\cdots (1-q^{2i_1+2i_2})}{(1-q^2)(1-q^4)\cdots (1-q^{2i_1})(1-q^2)(1-q^4)\cdots (1-q^{2i_2})} P_q (\Q [\mathcal{X}_{i_{3},1},\mathcal{X}_{i_{3},2}]).
\end{eqnarray*}
Since
\begin{equation*}
\frac{(1-q^2)(1-q^4)\cdots (1-q^{2i_1+2i_2})}{(1-q^2)(1-q^4)\cdots (1-q^{2i_1})(1-q^2)(1-q^4)\cdots (1-q^{2i_2})}=q^{i_1 i_2}\sum_{j=1}^{i_1 i_3 -i_1^2}p(j,i_3;i_1) \ q^{-i_1 i_3+i_1^2+2j},
\end{equation*}
the quotient ring $\Q [\mathcal{X}_{i_{3},1},\mathcal{X}_{i_{3},2},\mathcal{X}_{i_{1},3},\mathcal{X}_{i_{2},4}]/\left< x_{1,2}+X^{3;4}_{1,i_1;i_2},\cdots ,x_{i_3,2}+X^{3;4}_{i_3,i_1;i_2} \right> \left\{ -i_1 i_2 \right\}$ is isomorphic to 
\begin{equation*}
\bigoplus_{j=0}^{i_1 i_3-i_1^2}\left( \Q [\mathcal{X}_{i_{3},1},\mathcal{X}_{i_{3},2}] \left\{ -i_1 i_3+i_1^2 +2j\right\} \right)^{\oplus p(j,i_3;i_1)}
\end{equation*}
as a $\Z$-graded $\Q [\mathcal{X}_{i_{3},1},\mathcal{X}_{i_{3},2}]$-algebra.
Thus, we obtain the matrix factorization
\begin{equation*}
\bigoplus_{j=0}^{i_1 i_3 -i_1^2}\left( \c\acute{}\left( \input{fig-line-color-mf}\right)_n \{ -i_1 i_3 +i_1^2 +2j \}\right)^{\oplus p(j,i_3,i_1)}.
\end{equation*}

\noindent {\bf(2)} We have
\begin{eqnarray*}
&&\c\acute{}\left( \input{fig-bubble-color1-mf}\right)_n \\[-0.1em]
&=&\mathop{\boxtimes}_{j=1}^{i_3} K(V_{j,i_1;i_2}^{2;4,3};x_{j,3}-X_{j,i_1;i_2}^{2;4})_{\Q [\mathcal{X}_{i_{1},1},\mathcal{X}_{i_{1},2},\mathcal{X}_{i_{3},3},\mathcal{X}_{i_{2},4}]}\left\{ -i_1 i_2 \right\}
\mathop{\boxtimes}_{j=1}^{i_3} K(\Lambda_{j,i_1;i_2}^{1;4,3};X_{j,i_1;i_2}^{1;4}-x_{j,3})_{\Q[\mathcal{X}_{i_{1},1},\mathcal{X}_{i_{1},2},\mathcal{X}_{i_{3},3},\mathcal{X}_{i_{2},4}]}
\end{eqnarray*}

The potential of this matrix factorization does not include the variables of $\mathcal{X}_{i_3,3}$, $\mathcal{X}_{i_2,4}$ and the sequence 
\begin{equation*}
(X_{1,i_1;i_2}^{1;4}-x_{1,3},X_{2,i_1;i_2}^{1;4}-x_{2,3},\cdots ,X_{i_3,i_1;i_2}^{1;4}-x_{i_3,3})\left|_{(\mathcal{X}_{i_1,1},\mathcal{X}_{i_1,2})=
(\underline{0})}\right. 
=(x_{1,4}-x_{1,3},\cdots,x_{i_2,4}-x_{i_2,3},-x_{i_2+1,3},\cdots,-x_{i_3,3})
\end{equation*}
is regular in $\Q[\mathcal{X}_{i_{1},1},\mathcal{X}_{i_{1},2},\mathcal{X}_{i_{3},3},\mathcal{X}_{i_{2},4}]$. Thus we can apply Corollary \ref{cor2-11} 
to the matrix factorization. Then we have
\begin{eqnarray*}
&&\mathop{\boxtimes}_{j=1}^{i_3}K(\Lambda_{j,i_1;i_2}^{2;4,3};X_{j,i_1;i_2}^{1;4}-X_{j,i_1;i_2}^{2;4})_{\Q [\mathcal{X}_{i_{1},1},\mathcal{X}_{i_{1},2},\mathcal{X}_{i_{3},3},\mathcal{X}_{i_{2},4}]/\left< X_{1,i_1;i_2}^{1;4}-x_{1,3} ,\cdots ,X_{i_3,i_1;i_2}^{1;4}-x_{i_3,3} \right> }\left\{ -i_1 i_2 \right\}\\[-0.1em]
(\ast)&\simeq&\mathop{\boxtimes}_{j=1}^{i_3}K(\widetilde{\Lambda_{j,i_1;i_2}^{2;4,3}}\hspace{1mm};X_{j,i_1;i_2}^{1;4}-X_{j,i_1;i_2}^{2;4})_{\Q [\mathcal{X}_{i_{1},1},\mathcal{X}_{i_{1},2},\mathcal{X}_{i_{2},4}]}\left\{ -i_1 i_2 \right\} ,
\end{eqnarray*}
where 
\begin{equation*}
\widetilde{\Lambda_{j,i_1;i_2}^{2;4,3}}\hspace{1mm}=
\frac{F_{i_3}(X_{1,i_1;i_2}^{2;4},\cdots ,X_{j-1,i_1;i_2}^{2;4},X_{j,i_1;i_2}^{1;4},\cdots ,X_{i_3,i_1;i_2}^{1;4})-
F_{i_3}(X_{1,i_1;i_2}^{2;4},\cdots ,X_{j,i_1;i_2}^{2;4},X_{j+1,i_1;i_2}^{1;4},\cdots ,X_{i_3,i_1;i_2}^{1;4})}{X_{j,i_1;i_2}^{1;4}-X_{j,i_1;i_2}^{2;4}}.
\end{equation*}

Since $X_{j,i_1;i_2}^{1;4}-X_{j,i_1;i_2}^{2;4}$ is the $2j$-grading term of
\begin{equation*}
((x_{1,1}-x_{1,2})+(x_{2,1}-x_{2,2})+\cdots +(x_{i_1,1}-x_{i_1,2}))(1+x_{1,4}+x_{2,4}+\cdots +x_{i_2,4}),
\end{equation*}
the polynomials $X_{j,i_1;i_2}^{1;4}-X_{j,i_1;i_2}^{2;4}$ $(i_1+1\le j\le i_3)$ can be described 
as the linear sum of the polynomials $X_{j,i_1;i_2}^{1;4}-X_{j,i_1;i_2}^{2;4}$ $(1\le j\le i_1)$.
Using Proposition \ref{equiv} the matrix factorization $(\ast)$ is isomorphic to
\begin{eqnarray*}
&&\mathop{\boxtimes}_{j=1}^{i_1}K(\ast ;x_{j,1}-x_{j,2})_{\Q[\mathcal{X}_{i_1,1},\mathcal{X}_{i_1,2},\mathcal{X}_{i_{2},4}]}\boxtimes\\[-0.1em]
&&\mathop{\boxtimes}_{k=i_1+1}^{i_3}
(Q[\mathcal{X}_{i_{1},1},\mathcal{X}_{i_{1},2},\mathcal{X}_{i_{2},4}],Q[\mathcal{X}_{i_{1},1},\mathcal{X}_{i_{1},2},\mathcal{X}_{i_{2},4}]\{2k-n-1\},
\widetilde{\Lambda_{k,i_1;i_2}^{2;4,3}},0)\left\{ -i_1 i_2 \right\}.
\end{eqnarray*}
Since the potential of the partial matrix factorization 
$\displaystyle \mathop{\boxtimes}_{j=1}^{i_1}K(\ast ;x_{j,1}-x_{j,2})_{\Q[\mathcal{X}_{i_1,1},\mathcal{X}_{i_1,2},\mathcal{X}_{i_{2},4}]}$ equals 
$F_{i_1}(\mathcal{X}_{i_1,1})-F_{i_1}(\mathcal{X}_{i_1,2})$, 
by Theorem \ref{reg-eq} the partial matrix factorization is isomorphic to
\begin{equation*}
\c\acute{}\left( \input{fig-line-color1-mf}\right)_n\boxtimes (\Q[\mathcal{X}_{i_2,4}],0,0,0).
\end{equation*} 
The other partial matrix factorization $\mathop{\boxtimes}_{k=i_1+1}^{i_3}
(Q[\mathcal{X}_{i_{1},1},\mathcal{X}_{i_{1},2},\mathcal{X}_{i_{2},4}],Q[\mathcal{X}_{i_{1},1},\mathcal{X}_{i_{1},2},\mathcal{X}_{i_{2},4}]\{2k-n-1\},
\widetilde{\Lambda_{k,i_1;i_2}^{2;4,3}},0)\left\{ -i_1 i_2 \right\}$ is isomorphic to
\begin{equation*}
\mathop{\boxtimes}_{k=i_1+1}^{i_3}(0,Q[\mathcal{X}_{i_{1},1},\mathcal{X}_{i_{1},2},\mathcal{X}_{i_{2},4}]/\left< \widetilde{\Lambda_{k,i_1;i_2}^{2;4,3}}\hspace{1mm}\right>\{2k-n-1\},0,0).
\end{equation*}
Thus, we have
\begin{eqnarray*}
\c\acute{}\left( \input{fig-line-color1-mf}\right)_n\boxtimes
\mathop{\boxtimes}_{k=i_1+1}^{i_3}
(0,Q[\mathcal{X}_{i_{1},1},\mathcal{X}_{i_{1},2},\mathcal{X}_{i_{2},4}]/\left< \widetilde{\Lambda_{k,i_1;i_2}^{2;4,3}}\hspace{1mm}\right>\{2k-n-1\},0,0)
\left\{ -i_1 i_2 \right\}.\\
\end{eqnarray*}

Since we find that $\widetilde{\Lambda_{k,i_1;i_2}^{2;4,3}}
\left|_{(\mathcal{X}_{i_1,1},\mathcal{X}_{i_1,2})=(\underline{0})}\right.
=\frac{\partial F_{i_3}}{\partial x_{j}}(x_{1,4},\cdots,x_{i_1,4},0,\cdots,0)
\hspace{.5cm}
(i_1+1\le j\le i_3)$ form regular sequence,
these polynomials $\widetilde{\Lambda_{i_1+1,i_1;i_2}^{2;4,3}},\cdots,\widetilde{\Lambda_{i_3,i_1;i_2}^{2;4,3}} $ also form regular. The Poincar\'e series of the $\Z$-graded quotient ring 
\begin{equation*}
Q[\mathcal{X}_{i_{1},1},\mathcal{X}_{i_{1},2},\mathcal{X}_{i_{2},4}]/
\left< \widetilde{\Lambda_{i_1+1,i_1;i_2}^{2;4,3}}\hspace{1mm},\cdots ,\widetilde{\Lambda_{i_3,i_1;i_2}^{2;4,3}}\hspace{1mm}\right>\{\sum_{k=i_1+1}^{i_3}2k-n-1\}\left\{ -i_1 i_2 \right\}
\end{equation*}
equals to
\begin{eqnarray*}
&&\frac{(1-q^{2(n-i_1)})(1-q^{2(n-i_1-1)})\cdots (1-q^{2(n-i_3+1)})}{(1-q^2)(1-q^4)\cdots (1-q^{2i_2})} P_q (\Q [\mathcal{X}_{i_{1},1},\mathcal{X}_{i_{1},2}])
\prod_{k=i_1+1}^{i_3}q^{2k-n-1}\times q^{-i_1 i_2}\\[-0.1em]
&=&\left[n-i_1\atop i_2\right]P_q (\Q [\mathcal{X}_{i_{1},1},\mathcal{X}_{i_{1},2}]).
\end{eqnarray*}

Hence, the matrix factorization $\c\acute{}\left( \input{fig-bubble-color1-mf}\right)_n$ is isomorphic to\\[-1em]
\begin{equation*}
\bigoplus_{j=0}^{i_2(n-i_1)-i_2^2}\left( \c\acute{}\left( \input{fig-line-color1-mf}\right)_n \{ -i_2(n-i_1)+i_2^2 +2j \}\right)^{\oplus p(j,n-i_1,i_2)}\hspace{-1.5cm}\left< i_2 \right> .
\end{equation*}
\end{proof}
\begin{pro}\label{mat-equiv3}%proposition
There is the following isomorphism in $\HMF_{\Q[\mathcal{X}_{1,1},\mathcal{X}_{j,2},\mathcal{X}_{1,3},\mathcal{X}_{j,4}],F_1(\mathcal{X}_{1,1})+F_j(\mathcal{X}_{j,2})-F_1(\mathcal{X}_{1,3})-F_j(\mathcal{X}_{j,4})}$
\begin{eqnarray*}
{\bf(1)}\c\acute{}\left( \input{figsquare1j--1--2j-1--1--1j-mf}\right)_n &\simeq & \c\acute{}\left( \input{figsquare1j--j+1--1j-mf}\right)_n\bigoplus_{i=1}^{j-1} \c\acute{}\left( \input{figsquare1j-mf}\right)_n\{2i-j\}\\[-0.1em]
{\bf(2)}\c\acute{}\left( \input{figsquare1j--j+1--j1--j+1--1j-rev-mf}\right)_n &\simeq & \c\acute{}\left( \input{figsquare1j-rev-mf}\right)_n\bigoplus_{k=1}^{n-j-1} \c\acute{}\left( \input{figsquare1j--j-1--1j-rev-mf}\right)_n\{2k+j-n\} \left<1\right>.
\end{eqnarray*}
\end{pro}
\begin{proof}
{\bf(1)}We have
\begin{eqnarray*}
&&\c\acute{}\left( \input{figsquare1j--1--2j-1--1--1j-mf}\right)_n\\[-0.1em]
&=&
K\left(
\left( V^{1;6,5}_{1,1;1} \atop V^{1;6,5}_{2,1;1}\right);
\left(x_{1,1}+x_{1,6}-x_{1,5} \atop x_{1,1}x_{1,6}-x_{2,5}\right)
\right)_{\Q[\mathcal{X}_{1,1},\mathcal{X}_{1,6},\mathcal{X}_{2,5}]}\\[-0.1em]
&&\hspace{.5cm}\boxtimes 
K\left(
\left( \Lambda^{3;8,5}_{1,1;1} \atop \Lambda^{3;8,5}_{2,1;1}\right);
\left(x_{1,5}-x_{1,3}-x_{1,8} \atop x_{2,5}-x_{1,3}x_{1,8}\right)
\right)_{\Q[\mathcal{X}_{1,3},\mathcal{X}_{1,8},\mathcal{X}_{2,5}]}\hspace{-2cm}\{-1\}\\[-0.1em]
&&\hspace{1cm}\boxtimes
K\left(
\left( 
\begin{array}{c}
\Lambda^{6;7,2}_{1,1;j-1}\\[.5em]
\Lambda^{6;7,2}_{2,1;j-1}\\[.5em]
\Lambda^{6;7,2}_{3,1;j-1}\\[.5em]
\vdots\\[.5em]
\Lambda^{6;7,2}_{j-1,1;j-1}\\[.5em]
\Lambda^{6;7,2}_{j,1;j-1}
\end{array}
\right);
\left(
\begin{array}{c}
x_{1,2}-x_{1,6}-x_{1,7}\\[.6em]
x_{2,2}-x_{1,6}x_{1,7}-x_{2,7}\\[.6em]
x_{3,2}-x_{1,6}x_{2,7}-x_{3,7}\\[.6em]
\vdots\\[.6em]
x_{j-1,2}-x_{1,6}x_{j-2,7}-x_{j-1,7}\\[.6em]
x_{j,2}-x_{1,6}x_{j-1,7}
\end{array}
\right)
\right)_{\Q[\mathcal{X}_{j,2},\mathcal{X}_{1,6},\mathcal{X}_{j-1,7}]}\hspace{-2cm}\{-j+1\}\\[-0.1em]
&&\hspace{1.5cm}\boxtimes 
K\left(
\left( 
\begin{array}{c}
V^{8;7,4}_{1,1;j-1}\\[.5em]
V^{8;7,4}_{2,1;j-1}\\[.5em]
V^{8;7,4}_{3,1;j-1}\\[.5em]
\vdots\\[.5em]
V^{8;7,4}_{j-1,1;j-1}\\[.5em]
V^{8;7,4}_{j,1;j-1}
\end{array}
\right);
\left(
\begin{array}{c}
x_{1,7}+x_{1,8}-x_{1,4}\\[.6em]
x_{1,7}x_{1,8}+x_{2,7}-x_{2,4}\\[.6em]
x_{2,7}x_{1,8}+x_{3,7}-x_{3,4}\\[.6em]
\vdots\\[.6em]
x_{j-2,7}x_{1,8}+x_{j-1,7}-x_{j-1,4}\\[.6em]
x_{j-1,7}x_{1,8}-x_{j,4}
\end{array}
\right)
\right)_{\Q[\mathcal{X}_{1,8},\mathcal{X}_{j-1,7},\mathcal{X}_{j,2}].}
\end{eqnarray*}
We apply Corollary \ref{cor2-11} to the matrix factorization. Then we obtain
\begin{equation*}
K\left(
\left( 
\begin{array}{c}
\Lambda^{3;8,5}_{2,1;1}\\[.5em]
\Lambda^{6;7,2}_{j,1;j-1}\\[.5em]
V^{8;7,4}_{1,1;j-1}\\[.5em]
V^{8;7,4}_{2,1;j-1}\\[.5em]
V^{8;7,4}_{3,1;j-1}\\[.5em]
\vdots\\[.5em]
V^{8;7,4}_{j-1,1;j-1}\\[.5em]
V^{8;7,4}_{j,1;j-1}
\end{array}
\right);
\left(
\begin{array}{c}
x_{2,5}-x_{1,3}x_{1,8}\\[.6em]
x_{j,2}-x_{1,6}x_{j-1,7}\\[.6em]
x_{1,7}+x_{1,8}-x_{1,4}\\[.6em]
x_{1,7}x_{1,8}+x_{2,7}-x_{2,4}\\[.6em]
x_{2,7}x_{1,8}+x_{3,7}-x_{3,4}\\[.6em]
\vdots\\[.6em]
x_{j-2,7}x_{1,8}+x_{j-1,7}-x_{j-1,4}\\[.6em]
x_{j-1,7}x_{1,8}-x_{j,4}
\end{array}
\right)
\right)_{R_1}\{-j\},
\end{equation*}
where $R_1=\Q\left[
\begin{array}{c}
\mathcal{X}_{1,1},\mathcal{X}_{j,2},\mathcal{X}_{1,3},\mathcal{X}_{j,4},\\
\mathcal{X}_{2,5},\mathcal{X}_{1,6},\mathcal{X}_{j-1,7},\mathcal{X}_{1,8}
\end{array}
\right]\left/
\left< 
\begin{array}{c}
x_{1,1}+x_{1,6}-x_{1,5}, x_{1,1}x_{1,6}-x_{2,5}, x_{1,5}-x_{1,3}-x_{1,8}, \\
x_{1,2}-x_{1,6}-x_{1,7},  \cdots , x_{j,2}-x_{1,6}x_{j-1,7} 
\end{array}
\right>\right.$.\\
In the quotient ring, there are the following equations
\begin{eqnarray*}
x_{1,5}&=&x_{1,1}+x_{1,6},\\[-0.1em]
x_{2,5}&=&x_{1,1}x_{1,6},\\[-0.1em]
x_{1,8}&=&x_{1,1}-x_{1,3}+x_{1,6},\\[-0.1em]
x_{k,7}&=&\sum_{l=0}^{k}(-1)^lx_{1,6}^lx_{k-l,2}\,(=:A_k)\hspace{1cm}(k=1,2,\cdots ,j-1),
\end{eqnarray*}
where $x_{0,2}=1$.\\
Then we can consider 
$R_1\simeq \Q[\mathcal{X}_{1,1},\mathcal{X}_{j,2},\mathcal{X}_{1,3},\mathcal{X}_{j,4},\mathcal{X}_{1,6}]$. That is, the variables $x_{1,5}$, $x_{2,5}$, $x_{1,8}$ and $x_{k,7}$ can be removed from the quotient ring $R_1$ using the above equations. Therefore the matrix factorization is isomorphic to
\begin{equation*}
K\left(
\left( 
\begin{array}{c}
\widetilde{\Lambda^{3;8,5}_{2,1;1}}\\[.5em]
\widetilde{\Lambda^{6;7,2}_{j,1;j-1}}\\[.5em]
\widetilde{V^{8;7,4}_{1,1;j-1}}\\[.5em]
\widetilde{V^{8;7,4}_{2,1;j-1}}\\[.5em]
\vdots\\[.5em]
\widetilde{V^{8;7,4}_{j-1,1;j-1}}\\[.5em]
\widetilde{V^{8;7,4}_{j,1;j-1}}
\end{array}
\right);
\left(
\begin{array}{c}
(x_{1,6}-x_{1,3})(x_{1,1}-x_{1,3})\\[1em]
\sum_{l=0}^{j}(-1)^lx_{1,6}^l x_{j-l,2}\\[1em]
x_{1,1}+x_{1,2}-x_{1,3}-x_{1,4}\\[1em]
(x_{1,1}-x_{1,3})A_1+x_{2,2}-x_{2,4}\\[1em]
\vdots\\[1em]
(x_{1,1}-x_{1,3})A_{j-2}+x_{j-1,2}-x_{j-1,4}\\[1em]
(x_{1,1}-x_{1,3}+x_{1,6})A_{j-1}-x_{j,4}
\end{array}
\right)
\right)_{\Q[\mathcal{X}_{1,1},\mathcal{X}_{j,2},\mathcal{X}_{1,3},\mathcal{X}_{j,4},\mathcal{X}_{1,6}]}\hspace{-3cm}\{-j\}\hspace{3cm}.
\end{equation*}
Theorem \ref{exclude} is applied to $\sum_{l=0}^{j}(-1)^lx_{1,6}^l x_{j-l,2}$ of the matrix factorization.
Then we obtain
\begin{equation*}
K\left(
\left( 
\begin{array}{c}
\widehat{\Lambda^{3;8,5}_{2,1;1}}\\[.5em]
\widehat{V^{8;7,4}_{1,1;j-1}}\\[.5em]
\widehat{V^{8;7,4}_{2,1;j-1}}\\[.5em]
\vdots\\[.5em]
\widehat{V^{8;7,4}_{j,1;j-1}}
\end{array}
\right);
\left(
\begin{array}{c}
(x_{1,6}-x_{1,3})(x_{1,1}-x_{1,3})\\[.9em]
x_{1,1}+x_{1,2}-x_{1,3}-x_{1,4}\\[.9em]
(x_{1,1}-x_{1,3})A_1+x_{2,2}-x_{2,4}\\[.9em]
\vdots\\[.9em]
(x_{1,1}-x_{1,3})A_{j-1}+x_{j,2}-x_{j,4}
\end{array}
\right)
\right)_{R_1\acute{}}\{-j\},
\end{equation*}
where $R_1\acute{}=\Q[\mathcal{X}_{1,1},\mathcal{X}_{j,2},\mathcal{X}_{1,3},\mathcal{X}_{j,4},\mathcal{X}_{1,6}]\left/\left<\sum_{l=0}^{j}(-1)^lx_{1,6}^l x_{j-l,2}\right>\right.$, $\widehat{\Lambda^{3;8,5}_{2,1;1}}=\widetilde{\Lambda^{3;8,5}_{2,1;1}}$ and $\widehat{V^{8;7,4}_{i,1;j-1}}=\widetilde{V^{8;7,4}_{i,1;j-1}}$ in the quotient ring $R_1\acute{}$.
Since the polynomials $A_k$ are described as
\begin{eqnarray*}
A_k&=&\left(\sum_{l_1=0}^{k-1}(-1)^{k-l_1}x_{1,6}^{k-1-l_1}\left(\sum_{l_2=0}^{l_1}(-1)^{l_1-l_2}x_{1,3}^{l_1-l_2}x_{l_2,2}\right)\right)(x_{1,6}-x_{1,3})
+\sum_{l_3=0}^k (-1)^{k-l_3}x_{1,3}^{k-l_3}x_{l_3,2}\\
&=&u_{k-1}(x_{1,6}-x_1,3)+v_k,
\end{eqnarray*}
the above matrix factorization is isomorphic to
\begin{equation*}
K\left(
\left( 
\begin{array}{c}
\widehat{\Lambda^{3;8,5}_{2,1;1}}+\sum_{k=2}^{j}u_{k-2}\widehat{V^{8;7,4}_{k,1;j-1}}\\[.5em]
\widehat{V^{8;7,4}_{1,1;j-1}}\\[.5em]
\widehat{V^{8;7,4}_{2,1;j-1}}\\[.5em]
\vdots\\[.5em]
\widehat{V^{8;7,4}_{j,1;j-1}}
\end{array}
\right);
\left(
\begin{array}{c}
(x_{1,6}-x_{1,3})(x_{1,1}-x_{1,3})\\[.9em]
x_{1,1}+x_{1,2}-x_{1,3}-x_{1,4}\\[.9em]
(x_{1,1}-x_{1,3})v_1+x_{2,2}-x_{2,4}\\[.9em]
\vdots\\[.9em]
(x_{1,1}-x_{1,3})v_{j-1}+x_{j,2}-x_{j,4}
\end{array}
\right)
\right)_{R_1\acute{}}\{-j\}.
\end{equation*}

We can apply Corollary \ref{induce-sq1} to the matrix factorization, and then there are polynomials $a\in R_1\acute{}$ and $a_k\in \Q[\mathcal{X}_{1,1},\mathcal{X}_{j,2},\mathcal{X}_{1,3},\mathcal{X}_{j,4}](:=R_1\acute{}\,\acute{})$ ($k=1,\cdots j$) to give an isomorphism between the above matrix factorization and the following matrix factorization
\begin{equation}\label{total-mat1}
K\left(
\left( 
\begin{array}{c}
a\\[.1em]
a_1\\[.1em]
a_2\\[.1em]
\vdots\\[.1em]
a_j
\end{array}
\right);
\left(
\begin{array}{c}
(x_{1,6}-x_{1,3})(x_{1,1}-x_{1,3})\\[.1em]
x_{1,1}+x_{1,2}-x_{1,3}-x_{1,4}\\[.1em]
(x_{1,1}-x_{1,3})v_1+x_{2,2}-x_{2,4}\\[.1em]
\vdots\\[.1em]
(x_{1,1}-x_{1,3})v_{j-1}+x_{j,2}-x_{j,4}
\end{array}
\right)
\right)_{R_1\acute{}}\{-j\}.
\end{equation}
The partial matrix factorization $K(a;(x_{1,6}-x_{1,3})(x_{1,1}-x_{1,3}))$ has the potential 
\begin{equation*}
a(x_{1,6}-x_{1,3})(x_{1,1}-x_{1,3})=\omega-\sum_{i=1}^j a_i(-(x_{1,1}-x_{1,3})(\sum_{k=0}^{i-1} (-1)^{i-1-k}x_{1,3}^{i-1-k}x_{k,2})+x_{i,2}-x_{i,4})\in R_1\acute{}\,\acute{}.
\end{equation*}
Thus we have
\begin{eqnarray*}
&&k(a;(x_{1,6}-x_{1,3})(x_{1,1}-x_{1,3}))_{R_1\acute{}}\\[-0.1em]
\simeq&&
\left(
\hbox{
$
\left(
\begin{array}{r}
R_1\acute{}\,\acute{}\\
(x_{1,6}-x_{1,3})R_1\acute{}\,\acute{}\\
x_{1,6}(x_{1,6}-x_{1,3})R_1\acute{}\,\acute{}\\
\vdots\hspace{2cm}\\
x_{1,6}^{j-3}(x_{1,6}-x_{1,3})R_1\acute{}\,\acute{}\\
x_{1,6}^{j-2}(x_{1,6}-x_{1,3})R_1\acute{}\,\acute{}
\end{array}
\right)
$
}
,
\hbox{
$
\left(
\begin{array}{r}
R_1\acute{}\,\acute{}\{3-n\}\\
x_{1,6}R_1\acute{}\,\acute{}\{3-n\}\\
x_{1,6}^2R_1\acute{}\,\acute{}\{3-n\}\\
\vdots\hspace{2cm}\\
x_{1,6}^{j-2}R_1\acute{}\,\acute{}\{3-n\}\\
\alpha R_1\acute{}\,\acute{}\{3-n\}
\end{array}
\right)
$
}
,f_0,f_1
\right).
\end{eqnarray*}
where 
\begin{eqnarray*}
f_0&=&\left(
\begin{array}{cccc}
0&a(x_{1,6}-x_{1,3})&&\\
\vdots&&\ddots&\\
0&&&a(x_{1,6}-x_{1,3})\\
\displaystyle\frac{a}{\alpha}&0&\cdots&0
\end{array}
\right),\\
f_1&=&\left(
\begin{array}{cccc}
0&\cdots&0&\alpha (x_{1,6}-x_{1,3})(x_{1,1}-x_{1,3})\\[-0.1em]
x_{1,1}-x_{1,3}&&&0\\
&\ddots&&\vdots\\
&&x_{1,1}-x_{1,3}&0
\end{array}
\right),\\[-0.1em]
\alpha&=&\sum_{k_1=0}^{j-1}(-1)^{j-1-k_1}x_{1,6}^{j-1-k_1}\left(\sum_{k_2=0}^{k_1}(-1)^{k_1-k_2}x_{1,3}^{k_1-k_2}x_{k_2,2}\right).
\end{eqnarray*}
Remark that it is obvious to find
\begin{equation*}
\alpha (x_{1,6}-x_{1,3})=\sum_{k=0}^{j} (-1)^{j-k}x_{1,3}^{j-k}x_{k,2}.
\end{equation*}
Hence the partial matrix factorization splits into the following direct sum
\begin{eqnarray*}
&&\begin{array}{l}
(\xymatrix{(x_{1,6}-x_{1,3})R_1\acute{}\,\acute{}\ar[rr]^{a(x_{1,6}-x_{1,3})}&&R_1\acute{}\,\acute{}\{3-n\}\ar[rr]^{x_{1,1}-x_{1,3}}&&(x_{1,6}-x_{1,3})R_1\acute{}\,\acute{}})\\
\oplus
(\xymatrix{x_{1,6}(x_{1,6}-x_{1,3})R_1\acute{}\,\acute{}\ar[rr]^{a(x_{1,6}-x_{1,3})}&&x_{1,6}R_1\acute{}\,\acute{}\{3-n\}\ar[rr]^{x_{1,1}-x_{1,3}}&&x_{1,6}(x_{1,6}-x_{1,3})R_1\acute{}\,\acute{}})\\
\vdots\\
\oplus
(\xymatrix{x_{1,6}^{j-2}(x_{1,6}-x_{1,3})R_1\acute{}\,\acute{}\ar[rr]^{a(x_{1,6}-x_{1,3})}&&x_{1,6}^{j-2}R_1\acute{}\,\acute{}\{3-n\}\ar[rr]^{x_{1,1}-x_{1,3}}&&x_{1,6}^{j-2}(x_{1,6}-x_{1,3})R_1\acute{}\,\acute{}})\\
\oplus
(\xymatrix{R_1\acute{}\,\acute{}\ar[rrr]^(.4){\frac{a}{\alpha}}&&&\alpha R_1\acute{}\,\acute{}\{3-n\}\ar[rrr]^(.6){\alpha (x_{1,6}-x_{1,3})(x_{1,1}-x_{1,3})}&&&R_1\acute{}\,\acute{}}),
\end{array}\\[2em]
&\simeq&\begin{array}{l}
(\xymatrix{R_1\acute{}\,\acute{}\ar[rr]^{a(x_{1,6}-x_{1,3})}&&R_1\acute{}\,\acute{}\{1-n\}\ar[rr]^{x_{1,1}-x_{1,3}}&&R_1\acute{}\,\acute{}})\{2\}\\
\vdots\\
\oplus
(\xymatrix{R_1\acute{}\,\acute{}\ar[rr]^{a(x_{1,6}-x_{1,3})}&&R_1\acute{}\,\acute{}\{1-n\}\ar[rr]^{x_{1,1}-x_{1,3}}&&R_1\acute{}\,\acute{}})\{2j-2\}\\
\oplus
(\xymatrix{R_1\acute{}\,\acute{}\ar[rrr]^(.4){\frac{a}{\alpha}}&&& R_1\acute{}\,\acute{}\{2j+1-n\}\ar[rrr]^(.6){\alpha (x_{1,6}-x_{1,3})(x_{1,1}-x_{1,3})}&&&R_1\acute{}\,\acute{}}).
\end{array}\\[-0.1em]
&\simeq&
\bigoplus_{i=1}^{j-1}K(a(x_{1,6}-x_{1,3});x_{1,1}-x_{1,3})_{R_1\acute{}\,\acute{}}\{2i\}\oplus K(\frac{a}{\alpha};\alpha (x_{1,6}-x_{1,3})(x_{1,1}-x_{1,3}))_{R_1\acute{}\,\acute{}}
\end{eqnarray*}
By the decomposition and using Theorem \ref{reg-eq}, The total matrix factorization (\ref{total-mat1}) is isomorphic to
\begin{eqnarray*}
&&\hspace{-1cm}\bigoplus_{i=1}^{j-1}K\left(
\left( 
\begin{array}{c}
a(x_{1,6}-x_{1,3})\\[.1em]
a_1\\[.1em]
a_2\\[.1em]
\vdots\\[.1em]
a_j
\end{array}
\right);
\left(
\begin{array}{c}
x_{1,1}-x_{1,3}\\[.1em]
x_{1,1}+x_{1,2}-x_{1,3}-x_{1,4}\\[.1em]
(x_{1,1}-x_{1,3})(\sum_{k=0}^{1} (-1)^{1-k}x_{1,3}^{1-k}x_{k,2})+x_{2,2}-x_{2,4}\\[.1em]
\vdots\\[.1em]
(x_{1,1}-x_{1,3})(\sum_{k=0}^{j-1} (-1)^{j-1-k}x_{1,3}^{j-1-k}x_{k,2})+x_{j,2}-x_{j,4}
\end{array}
\right)
\right)_{R_1\acute{}\,\acute{}}\{2i-j\}\\[-0.1em]
&&\bigoplus
K\left(
\left( 
\begin{array}{c}
a_1\\[.1em]
a_2\\[.1em]
\vdots\\[.1em]
a_j\\[.1em]
\frac{a}{\alpha}
\end{array}
\right);
\left(
\begin{array}{c}
x_{1,1}+x_{1,2}-x_{1,3}-x_{1,4}\\[.1em]
(x_{1,1}-x_{1,3})(\sum_{k=0}^{1} (-1)^{1-k}x_{1,3}^{1-k}x_{k,2})+x_{2,2}-x_{2,4}\\[.1em]
\vdots\\[.1em]
(x_{1,1}-x_{1,3})(\sum_{k=0}^{j-1} (-1)^{j-1-k}x_{1,3}^{j-1-k}x_{k,2})+x_{j,2}-x_{j,4}\\[.1em]
(x_{1,1}-x_{1,3})(\sum_{k=0}^{j} (-1)^{j-k}x_{1,3}^{j-k}x_{k,2})
\end{array}
\right)
\right)_{R_1\acute{}\,\acute{}}\{-j\}\\[-0.1em]
&&\hspace{-1cm}\simeq
\bigoplus_{i=1}^{j-1}K\left(
\left( 
\begin{array}{c}
a(x_{1,6}-x_{1,3})+\sum_{l=1}^{j}(\sum_{k=0}^{l-1} (-1)^{l-2-k}x_{1,3}^{l-1-k}x_{k,2})a_l\\[.1em]
a_1\\[.1em]
a_2\\[.1em]
\vdots\\[.1em]
a_j
\end{array}
\right);
\left(
\begin{array}{c}
x_{1,1}-x_{1,3}\\[.1em]
x_{1,2}-x_{1,4}\\[.1em]
x_{2,2}-x_{2,4}\\[.1em]
\vdots\\[.1em]
x_{j,2}-x_{j,4}
\end{array}
\right)
\right)_{R_1\acute{}\,\acute{}}\{2i-j\}\\[-0.1em]
&&\bigoplus
K\left(
\left( 
\begin{array}{c}
\sum_{k=1}^{j+1}(-1)^{k-1}a_k x_{1,3}^{k-1}\\[.1em]
\sum_{k=2}^{j+1}(-1)^{k-2}a_k x_{1,3}^{k-2}\\[.1em]
\vdots\\[.1em]
a_{j-1}-x_{1,3}a_{j}+x_{1,3}^2a_{j+1}\\[.1em]
a_j-x_{1,3}a_{j+1}\\[.1em]
a_{j+1}
\end{array}
\right);
\left(
\begin{array}{c}
x_{1,1}+x_{1,2}-x_{1,3}-x_{1,4}\\[.1em]
x_{1,1}x_{1,2}-x_{1,3}x_{1,4}+x_{2,2}-x_{2,4}\\[.1em]
\vdots\\[.1em]
x_{1,1}x_{j-2,2}-x_{1,3}x_{j-2,4}+x_{j-1,2}-x_{j-1,4}\\[.1em]
x_{1,1}x_{j-1,2}-x_{1,3}x_{j-1,4}+x_{j,2}-x_{j,4}\\[.1em]
x_{1,1}x_{j,2}-x_{1,3}x_{j,4}
\end{array}
\right)
\right)_{R_1\acute{}\,\acute{}}\{-j\}
\end{eqnarray*}
where $a_{j+1}=\frac{a}{\alpha}.$\\
Using Theorem \ref{reg-eq}, we find the above matrix factorization is isomorphic to
\begin{equation*}
\c\acute{}\left( \input{figsquare1j--j+1--1j-mf}\right)\bigoplus_{i=1}^{j-1} \c\acute{}\left( \input{figsquare1j-mf}\right)_n\{2i-j\}
\end{equation*}

%proof of (2)
{\bf(2)}We have
\begin{eqnarray*}
&&
K\left(\left(
\begin{array}{c}
\Lambda_{1,1;j}^{1;5,6}\\[.5em]
\Lambda_{2,1;j}^{1;5,6}\\[.5em]
\vdots\\[.5em]
\Lambda_{j,1;j}^{1;5,6}\\[.5em]
\Lambda_{j+1,1;j}^{1;5,6}
\end{array}
\right)
;
\left(
\begin{array}{c}
x_{1,6}-x_{1,1}-x_{1,5}\\[.6em]
x_{2,6}-x_{1,1}x_{1,5}-x_{2,5}\\[.6em]
\vdots\\[.6em]
x_{j,6}-x_{1,1}x_{j-1,5}-x_{j,5}\\[.6em]
x_{j+1,6}-x_{1,1}x_{j,5}-x_{j+1,5}
\end{array}
\right)
\right)_{\Q[\mathcal{X}_{1,1}\mathcal{X}_{j,5}\mathcal{X}_{j+1,6}]}\hspace{-2cm}\{-j\}\\
&&\hspace{.5cm}\boxtimes
K\left(\left(
\begin{array}{c}
V_{1,1;j}^{3;5,8}\\[.5em]
V_{2,1;j}^{3;5,8}\\[.5em]
\vdots\\[.5em]
V_{j,1;j}^{3;5,8}\\[.5em]
V_{j+1,1;j}^{3;5,8}
\end{array}
\right)
;
\left(
\begin{array}{c}
x_{1,3}+x_{1,5}-x_{1,8}\\[.6em]
x_{1,3}x_{1,5}+x_{2,5}-x_{2,8}\\[.6em]
\vdots\\[.6em]
x_{1,3}x_{j-1,5}+x_{j,5}-x_{j,8}\\[.6em]
x_{1,3}x_{j,5}-x_{j+1,8}
\end{array}
\right)
\right)_{\Q[\mathcal{X}_{1,3}\mathcal{X}_{j,5}\mathcal{X}_{j+1,8}]}\\
&&\hspace{1cm}\boxtimes
K\left(\left(
\begin{array}{c}
V_{1,1;j}^{7;2,6}\\[.5em]
V_{2,1;j}^{7;2,6}\\[.5em]
\vdots\\[.5em]
V_{j,1;j}^{7;2,6}\\[.5em]
V_{j+1,1;j}^{7;2,6}
\end{array}
\right)
;
\left(
\begin{array}{c}
x_{1,7}+x_{1,2}-x_{1,6}\\[.6em]
x_{1,7}x_{1,2}+x_{2,2}-x_{2,6}\\[.6em]
\vdots\\[.6em]
x_{1,7}x_{j-1,2}+x_{j,2}-x_{j,6}\\[.6em]
x_{1,7}x_{j,2}-x_{j+1,6}
\end{array}
\right)
\right)_{\Q[\mathcal{X}_{1,7}\mathcal{X}_{j,2}\mathcal{X}_{j+1,6}]}\\
&&\hspace{1.5cm}\boxtimes
K\left(\left(
\begin{array}{c}
\Lambda_{1,1;j}^{7;4,8}\\[.5em]
\Lambda_{2,1;j}^{7;4,8}\\[.5em]
\vdots\\[.5em]
\Lambda_{j,1;j}^{7;4,8}\\[.5em]
\Lambda_{j+1,1;j}^{7;4,8}
\end{array}
\right)
;
\left(
\begin{array}{c}
x_{1,8}-x_{1,7}-x_{1,4}\\[.6em]
x_{2,8}-x_{1,7}x_{1,4}-x_{2,4}\\[.6em]
\vdots\\[.6em]
x_{j,8}-x_{1,7}x_{j-1,4}-x_{j,4}\\[.6em]
x_{j+1,8}-x_{1,7}x_{j,4}
\end{array}
\right)
\right)_{\Q[\mathcal{X}_{1,7}\mathcal{X}_{j,4}\mathcal{X}_{j+1,8}]}\hspace{-2cm}\{-j\}.
\end{eqnarray*}
We apply Corollary \ref{cor2-11} to the matrix factorization. Then we obtain
\begin{equation*}
K\left(\left(
\begin{array}{c}
\Lambda_{1,1;j}^{1;5,6}\\[.5em]
\vdots\\[.5em]
\Lambda_{j+1,1;j}^{1;5,6}\\[.5em]
V_{j+1,1;j}^{3;5,8}
\end{array}
\right)
;
\left(
\begin{array}{c}
x_{1,6}-x_{1,1}-x_{1,5}\\[.6em]
\vdots\\[1em]
x_{j+1,6}-x_{1,1}x_{j,5}\\[.6em]
x_{1,3}x_{j,5}-x_{j+1,8}
\end{array}
\right)
\right)_{R_2}\{-2j\},
\end{equation*}
where 
$
R_2=\Q\left[
\begin{array}{c}
\mathcal{X}_{1,1},\mathcal{X}_{j,2},\mathcal{X}_{1,3},\mathcal{X}_{j,4},\\
\mathcal{X}_{j,5},\mathcal{X}_{j+1,6},\mathcal{X}_{1,7},\mathcal{X}_{j+1,8}
\end{array}
\right]
\left/\left<
\begin{array}{c}
x_{1,3}+x_{1,5}-x_{1,8},
\cdots,
x_{1,3}x_{j-1,5}+x_{j,5}-x_{j,8},\\
x_{1,7}+x_{1,2}-x_{1,6},
\cdots,
x_{1,7}x_{j,2}-x_{j+1,6},\\
x_{1,8}-x_{1,7}-x_{1,4},
\cdots,
x_{j+1,8}-x_{1,7}x_{j,4}
\end{array}
\right>\right.
$\\

In the quotient ring, there are equations
\begin{eqnarray*}
x_{k,5}&=&\left(\sum_{l=0}^{k-1}(-1)^{k-1-l}x_{1,3}^{k-1-l}x_{l,4}\right)x_{1,7}+\sum_{l=0}^{k}(-1)^{k-l}x_{1,3}^{k-l}x_{l,4}(=:B_k)\hspace{1cm}(k=1,\cdots,j),\\[-0.1em]
x_{1,6}&=&x_{1,7}+x_{1,2},\\[-0.1em]
&&\vdots\\[-0.1em]
x_{j+1,6}&=&x_{1,7}x_{j,2},\\[-0.1em]
x_{1,8}&=&x_{1,7}+x_{1,4},\\[-0.1em]
&&\vdots\\[-0.1em]
x_{j+1,8}&=&x_{1,7}x_{j,4}.
\end{eqnarray*}
Using the above equations, we find the matrix factorization is isomorphic to
\begin{eqnarray*}
&&
K\left(\left(
\begin{array}{c}
\widetilde{\Lambda_{1,1;j}^{1;5,6}}\\[.5em]
\widetilde{\Lambda_{2,1;j}^{1;5,6}}\\[.5em]
\vdots\\[.5em]
\widetilde{\Lambda_{j,1;j}^{1;5,6}}\\[.5em]
\widetilde{\Lambda_{j+1,1;j}^{1;5,6}}\\[.5em]
\widetilde{V_{j+1,1;j}^{3;5,8}}
\end{array}
\right)
;
\left(
\begin{array}{c}
s_1\\[1em]
s_1 x_{1,7}+s_2\\[1em]
\vdots\\[1em]
s_{j-1} x_{1,7}+s_j\\[1em]
x_{1,7}x_{j,2}-x_{1,1}B_j\\[1em]
x_{1,3}B_j-x_{1,7}x_{j,4}
\end{array}
\right)
\right)_{\Q[\mathcal{X}_{1,1},\mathcal{X}_{j,2},\mathcal{X}_{1,3},\mathcal{X}_{j,4},\mathcal{X}_{1,7}]}\hspace{-2.5cm}\{-2j\},\\
&\simeq&
K\left(\left(
\begin{array}{c}
\widetilde{\Lambda_{1,1;j}^{1;5,6}}+x_{1,7}\widetilde{\Lambda_{2,1;j}^{1;5,6}}\\[.5em]
\vdots\\[.5em]
\widetilde{\Lambda_{j-1,1;j}^{1;5,6}}+x_{1,7}\widetilde{\Lambda_{j,1;j}^{1;5,6}}\\[.5em]
\widetilde{\Lambda_{j,1;j}^{1;5,6}}\\[.5em]
\widetilde{\Lambda_{j+1,1;j}^{1;5,6}}\\[.5em]
\widetilde{V_{j+1,1;j}^{3;5,8}}
\end{array}
\right)
;
\left(
\begin{array}{c}
s_1\\[1em]
\vdots\\[1em]
s_{j-1}\\[1em]
s_j\\[1em]
x_{1,7}x_{j,2}-x_{1,1}B_j\\[1em]
(x_{1,3}-x_{1,7})(\sum_{l=0}^{j}(-1)^{j-l}x_{1,3}^{j-l}x_{l,4})
\end{array}
\right)
\right)_{\Q[\mathcal{X}_{1,1},\mathcal{X}_{j,2},\mathcal{X}_{1,3},\mathcal{X}_{j,4},\mathcal{X}_{1,7}]}\hspace{-2.5cm}\{-2j\},
\end{eqnarray*}
where $s_k=x_{k,2}+\sum_{l=0}^{k}(-1)^{k+1-l}x_{1,3}^{k-l}X_{l,1;j}^{1;4}$ and $\widetilde{\Lambda_{k,1;j}^{1;5,6}},$ $\widetilde{V_{j+1,1;j}^{3;5,8}}$ are derived from $\Lambda_{k,1;j}^{1;5,6},$ $V_{j+1,1;j}^{3;5,8}$ using the equations.\\
The polynomial $\widetilde{\Lambda_{j+1,1;j}^{1;5,6}}$ can be described by
\begin{eqnarray*}
\widetilde{\Lambda_{j+1,1;j}^{1;5,6}}&=&\frac{F_{j+1}(X_{1,1;j}^{1;5},\cdots,X_{j,1;j}^{1;5},x_{j+1,6})-F_{j+1}(X_{1,1;j}^{1;5},\cdots,X_{j,1;j}^{1;5},X_{j+1,1;j}^{1;5})}{x_{j+1,6}-X_{j+1,1;j}^{1;5}}\\[-0.1em]
&=&c_0\left(X_{1,1;j}^{1;5}\right)^{n-j}+\cdots\\[-0.1em]
&=&c_0(x_{1,7}+x_{1,1}-x_{1,3}+x_{1,4})^{n-j}+\cdots\\[-0.1em]
&=&c_0 x_{1,7}^{n-j}+c_1 x_{1,7}^{n-j-1} + \cdots +c_{n-j},
\end{eqnarray*}
where $c_0 \in \Q$ and $c_k \in\Q[\mathcal{X}_{1,1},\mathcal{X}_{j,2},\mathcal{X}_{1,3},\mathcal{X}_{j,4}](:=R_2\acute{}\,\acute{})$ ($k=1,\cdots,n-j$). Using Theorem \ref{exclude}, we have
\begin{eqnarray*}
&&
K\left(\left(
\begin{array}{c}
\widetilde{\Lambda_{1,1;j}^{1;5,6}}+x_{1,7}\widetilde{\Lambda_{2,1;j}^{1;5,6}}\\[.5em]
\vdots\\[.5em]
\widetilde{\Lambda_{j-1,1;j}^{1;5,6}}+x_{1,7}\widetilde{\Lambda_{j,1;j}^{1;5,6}}\\[.5em]
\widetilde{\Lambda_{j,1;j}^{1;5,6}}\\[.5em]
x_{1,7}x_{j,2}-x_{1,1}B_j\\[1em]
\widetilde{V_{j+1,1;j}^{3;5,8}}
\end{array}
\right)
;
\left(
\begin{array}{c}
s_1\\[1em]
\vdots\\[1em]
s_{j-1}\\[1em]
s_j\\[1em]
\widetilde{\Lambda_{j+1,1;j}^{1;5,6}}\\[.5em]
(x_{1,3}-x_{1,7})(\sum_{l=0}^{j}(-1)^{j-l}x_{1,3}^{j-l}x_{l,4})
\end{array}
\right)
\right)_{\Q[\mathcal{X}_{1,1},\mathcal{X}_{j,2},\mathcal{X}_{1,3},\mathcal{X}_{j,4},\mathcal{X}_{1,7}]}\hspace{-2.5cm}\{-2j\}\{2j+1-n\}\left<1\right>,\\[-0.1em]
&\simeq&
K\left(\left(
\begin{array}{c}
\widetilde{\Lambda_{1,1;j}^{1;5,6}}+x_{1,7}\widetilde{\Lambda_{2,1;j}^{1;5,6}}\\[.5em]
\vdots\\[.5em]
\widetilde{\Lambda_{j-1,1;j}^{1;5,6}}+x_{1,7}\widetilde{\Lambda_{j,1;j}^{1;5,6}}\\[.5em]
\widetilde{\Lambda_{j,1;j}^{1;5,6}}\\[.5em]
\widetilde{V_{j+1,1;j}^{3;5,8}}
\end{array}
\right)
;
\left(
\begin{array}{c}
s_1\\[1em]
\vdots\\[1em]
s_{j-1}\\[1em]
s_j\\[1em]
(x_{1,3}-x_{1,7})(\sum_{l=0}^{j}(-1)^{j-l}x_{1,3}^{j-l}x_{l,4})
\end{array}
\right)
\right)_{R_2\acute{}}\{1-n\}\left<1\right>,
\end{eqnarray*}
where $R_2\acute{}=\Q[\mathcal{X}_{1,1},\mathcal{X}_{j,2},\mathcal{X}_{1,3},\mathcal{X}_{j,4},\mathcal{X}_{1,7}]/\widetilde{\Lambda_{j+1,1;j}^{1;5,6}}$. Using Corollary \ref{induce-sq1}, there are polynomials $b_k \in R_2\acute{}\,\acute{} $ ($k=1,\cdots,j$) and $b\in R_2\acute{}$ to give the following matrix factorization which is isomorphic to the above matrix factorization
\begin{equation*}
K\left(\left(
\begin{array}{c}
b_1\\[.5em]
\vdots\\[.5em]
b_j\\[.6em]
b
\end{array}
\right)
;
\left(
\begin{array}{c}
s_1\\[.5em]
\vdots\\[.5em]
s_j\\[.5em]
(x_{1,3}-x_{1,7})(\sum_{l=0}^{j}(-1)^{j-l}x_{1,3}^{j-l}x_{l,4})
\end{array}
\right)
\right)_{R_2\acute{}}\{1-n\}\left<1\right>.
\end{equation*}
The partial matrix factorization $K(b;(x_{1,3}-x_{1,7})(\sum_{l=0}^{j}(-1)^{j-l}x_{1,3}^{j-l}x_{l,4}))_{R_2\acute{}}$ is described as
\begin{eqnarray*}
&&K\left(b;(x_{1,3}-x_{1,7})\left(\sum_{l=0}^{j}(-1)^{j-l}x_{1,3}^{j-l}x_{l,4}\right)\right)_{R_2\acute{}}\\
&\simeq&
\left(
\hbox{
$
\left(
\begin{array}{r}
R_2\acute{}\,\acute{}\\
(x_{1,3}-x_{1,7})R_2\acute{}\,\acute{}\\
x_{1,7}(x_{1,3}-x_{1,7})R_2\acute{}\,\acute{}\\
\vdots\hspace{2cm}\\
x_{1,7}^{n-j-3}(x_{1,3}-x_{1,7})R_2\acute{}\,\acute{}\\
x_{1,7}^{n-j-2}(x_{1,3}-x_{1,7})R_2\acute{}\,\acute{}
\end{array}
\right)
$
}
,
\hbox{
$
\left(
\begin{array}{r}
R_2\acute{}\,\acute{}\{3-n\}\\
x_{1,7} R_2\acute{}\,\acute{}\{3-n\}\\
x_{1,7}^2 R_2\acute{}\,\acute{}\{3-n\}\\
\vdots\hspace{2cm}\\
x_{1,7}^{n-j-2} R_2\acute{}\,\acute{}\{3-n\}\\
\beta R_2\acute{}\,\acute{}\{3-n\}
\end{array}
\right)
$
}
,g_0,g_1\right),
\end{eqnarray*}
where 
\begin{eqnarray*}
g_0&=&\left(
\begin{array}{cccc}
0&b(x_{1,3}-x_{1,7})&&\\
\vdots&&\ddots&\\
0&&&b(x_{1,3}-x_{1,7})\\
\frac{b}{\beta}&0&\cdots&0
\end{array}
\right),\\[-0.1em]
g_1&=&\left(
\begin{array}{cccc}
0&\cdots&0&\displaystyle\beta (x_{1,3}-x_{1,7})(\sum_{l=0}^{j}(-1)^{j-l}x_{1,3}^{j-l}x_{l,4})\\
\displaystyle\sum_{l=0}^{j}(-1)^{j-l}x_{1,3}^{j-l}x_{l,4}&&&0\\
&\ddots&&\vdots\\
&&\displaystyle\sum_{l=0}^{j}(-1)^{j-l}x_{1,3}^{j-l}x_{l,4}&0
\end{array}
\right),\\[-0.1em]
\beta&=&\sum_{k=0}^{n-j-1}x_{1,7}^{n-j-1-k}\left(\sum_{l=0}^{k}x_{1,3}^{k-l}c_l\right),\\
R_2\acute{}\,\acute{}&=&\Q[\mathcal{X}_{1,1},\mathcal{X}_{j,2},\mathcal{X}_{1,3},\mathcal{X}_{j,4}].
\end{eqnarray*}
Remark that we have
\begin{eqnarray*}
\beta(x_{1,3}-x_{1,7})&=&c_0x_{1,3}^{n-j}+c_1x_{1,3}^{n-j-1}+\cdots +c_{n-j},\\[-0.1em]
&\simeq&\widetilde{\Lambda_{j+1,1;j}^{1;5,6}}\left|_{x_{1,7}\to x_{1,3}}\right. .
\end{eqnarray*}
Therefore the partial matrix factorization is isomorphic to
\begin{eqnarray*}
&&\bigoplus_{k=1}^{n-j-1}K\left(b(x_{1,3}-x_{1,7});\sum_{l=0}^{j}(-1)^{j-l}x_{1,3}^{j-l}x_{l,4}\right)_{R_2\acute{}\,\acute{}}\{2k\}\\[-0.1em]
&&\hspace{1cm}\bigoplus K\left(\frac{b}{\beta};\beta (x_{1,3}-x_{1,7})\left(\sum_{l=0}^{j}(-1)^{j-l}x_{1,3}^{j-l}x_{l,4}\right)\right)_{R_2\acute{}\,\acute{}}.
\end{eqnarray*}
Then the total matrix factorization is isomorphic to
\begin{eqnarray}\label{total-mat2}
&&\bigoplus_{k=1}^{n-j-1}
K\left(\left(
\begin{array}{c}
b_1\\[.5em]
\vdots\\[.5em]
b_j\\[.6em]
b(x_{1,3}-x_{1,7})
\end{array}
\right)
;
\left(
\begin{array}{c}
s_1\\[.5em]
\vdots\\[.5em]
s_j\\[.5em]
\sum_{l=0}^{j}(-1)^{j-l}x_{1,3}^{j-l}x_{l,4}
\end{array}
\right)
\right)_{R_2\acute{},\acute{}}\{2k+1-n\}\left<1\right>\\[-0.1em]
\label{total-mat3}&&\hspace{1cm}\bigoplus
K\left(\left(
\begin{array}{c}
b_1\\[.5em]
\vdots\\[.5em]
b_j\\[.6em]
\frac{b}{\beta}
\end{array}
\right)
;
\left(
\begin{array}{c}
s_1\\[.5em]
\vdots\\[.5em]
s_j\\[.5em]
\widetilde{\Lambda_{j+1,1;j}^{1;5,6}}\left|_{x_{1,7}\to x_{1,3}}\right.\left(\sum_{l=0}^{j}(-1)^{j-l}x_{1,3}^{j-l}x_{l,4}\right)
\end{array}
\right)
\right)_{R_2\acute{},\acute{}}\{1-n\}\left<1\right>.
\end{eqnarray}
On the other hand, we have
\begin{eqnarray*}
\c\acute{}\left(\input{figsquare1j--j-1--1j-rev-mf}\right)_n&=&
K\left(\left(
\begin{array}{c}
\Lambda_{1,1:j-1}^{1;9,2}\\[.5em]
\vdots\\[.5em]
\Lambda_{j,1:j-1}^{1;9,2}
\end{array}
\right)
;
\left(
\begin{array}{c}
x_{1,2}-X_{1,1;j-1}^{1;9}\\[.5em]
\vdots\\[.5em]
x_{j,2}-X_{j,1;j-1}^{1;9}
\end{array}
\right)
\right)_{\Q[\mathcal{X}_{1,1},\mathcal{X}_{j,2},\mathcal{X}_{j-1,9}]}\hspace{-2cm}\{1-j\}\\
&&\hspace{1cm}\boxtimes
K\left(\left(
\begin{array}{c}
V_{1,1:j-1}^{3;9,4}\\[.5em]
\vdots\\[.5em]
V_{j,1:j-1}^{3;9,4}
\end{array}
\right)
;
\left(
\begin{array}{c}
X_{1,1;j-1}^{3;9}-x_{1,4}\\[.5em]
\vdots\\[.5em]
X_{j,1;j-1}^{3;9}-x_{j,4}
\end{array}
\right)
\right)_{\Q[\mathcal{X}_{1,3},\mathcal{X}_{j,4},\mathcal{X}_{j-1,9}]}\\
&\simeq&
K\left(\left(
\begin{array}{c}
\widetilde{\Lambda_{1,1:j-1}^{1;9,2}}\\[.5em]
\vdots\\[.5em]
\widetilde{\Lambda_{j-2,1:j-1}^{1;9,2}}\\[.5em]
\widetilde{\Lambda_{j-1,1:j-1}^{1;9,2}}\\[.5em]
\widetilde{\Lambda_{j,1:j-1}^{1;9,2}}
\end{array}
\right)
;
\left(
\begin{array}{c}
s_1\\[1em]
\vdots\\[1em]
s_{j-1}\\[1em]
s_j+\sum_{l=0}^{j}(-1)^{j-l}x_{1,3}^{j-l}x_{l,4}\\[1em]
-\sum_{l=0}^{j}(-1)^{j-l}x_{1,3}^{j-l}x_{l,4}
\end{array}
\right)
\right)_{R_2\acute{}\,\acute{}}\{1-j\}\\
&\simeq&
K\left(\left(
\begin{array}{c}
\widetilde{\Lambda_{1,1:j-1}^{1;9,2}}\\[.5em]
\vdots\\[.5em]
\widetilde{\Lambda_{j-1,1:j-1}^{1;9,2}}\\[.5em]
\widetilde{\Lambda_{j,1:j-1}^{1;9,2}}-\widetilde{\Lambda_{j-1,1:j-1}^{1;9,2}}
\end{array}
\right)
;
\left(
\begin{array}{c}
s_1\\[1em]
\vdots\\[1em]
s_j\\[1em]
-\sum_{l=0}^{j}(-1)^{j-l}x_{1,3}^{j-l}x_{l,4}
\end{array}
\right)
\right)_{R_2\acute{}\,\acute{}}\{1-j\}
.
\end{eqnarray*}
Thus using Theorem \ref{reg-eq}, the partial matrix factorization (\ref{total-mat2}) is isomorphic to
\begin{equation*}
\bigoplus_{k=1}^{n-j-1}
\c\acute{}\left(\input{figsquare1j--j-1--1j-rev-mf}\right)_n
\{2k+j-n\}\left<1\right>.
\end{equation*}
Moreover we have
\begin{eqnarray*}
&&\widetilde{\Lambda_{j+1,1;j}^{1;5,6}}\left|_{x_{1,7}\to x_{1,3}}\right.\left(\sum_{l=0}^{j}(-1)^{j-l}x_{1,3}^{j-l}x_{l,4}\right)\\
&=&\left.\frac{F_{j+1}(X_{1,1;j}^{1;5},\cdots,X_{j,1;j}^{1;5},x_{j+1,6})-F_{j+1}(X_{1,1;j}^{1;5},\cdots,X_{j,1;j}^{1;5},X_{j+1,1;j}^{1;5})}{x_{j+1,6}-X_{j+1,1;j}^{1;5}}\right|_{x_{1,7}\to x_{1,3}}\left(\sum_{l=0}^{j}(-1)^{j-l}x_{1,3}^{j-l}x_{l,4}\right).
\end{eqnarray*}
Since we have
\begin{eqnarray*}
\left.X_{k,1;j}^{1;5}\right|_{x_{1,7}\to x_{1,3}}&=&\left.x_{1,1}x_{k-1,5}+x_{k,5}\right|_{x_{1,7}\to x_{1,3}}\\
&=&x_{1,1}\left(\left(\sum_{l=0}^{k-2}(-1)^{k-2-l}x_{1,3}^{k-2-l}x_{l,4}\right)x_{1,7}+\sum_{l=0}^{k-1}(-1)^{k-1-l}x_{1,3}^{k-1-l}x_{l,4}\right)\\
&&+\left.\left(\sum_{l=0}^{k-1}(-1)^{k-1-l}x_{1,3}^{k-1-l}x_{l,4}\right)x_{1,7}+\sum_{l=0}^{k}(-1)^{k-l}x_{1,3}^{k-l}x_{l,4}\right|_{x_{1,7}\to x_{1,3}}\\
&=&x_{1,1}x_{k-1,4}+x_{k,4}=X_{k,1;j}^{1;4},
\end{eqnarray*}
$\widetilde{\Lambda_{j+1,1;j}^{1;5,6}}\left|_{x_{1,7}\to x_{1,3}}\right.\left(\sum_{l=0}^{j}(-1)^{j+1-l}x_{1,3}^{j-l}x_{l,4}\right)$ equals
\begin{eqnarray*}
&&\frac{F_{j+1}(X_{1,1;j}^{1;4},\cdots,X_{j,1;j}^{1;4},x_{1,3}x_{j,2})-F_{j+1}(X_{1,1;j}^{1;4},\cdots,X_{j,1;j}^{1;4},X_{j+1,1;j}^{1;4})}{x_{j+1,6}-X_{j+1,1;j}^{1;4}}\left(\sum_{l=0}^{j}(-1)^{j-l}x_{1,3}^{j-l}x_{l,4}\right)\\
&\equiv&\frac{F_{j+1}(X_{1,1;j}^{1;4},\cdots,X_{j,1;j}^{1;4},x_{1,3}(\sum_{l=0}^{j}(-1)^{j+1-l}x_{1,3}^{j-l}X_{l,1;j}^{1;4}))-F_{j+1}(X_{1,1;j}^{1;4},\cdots,X_{j,1;j}^{1;4},X_{j+1,1;j}^{1;4})}{x_{1,3}(\sum_{l=0}^{j}(-1)^{j-l}x_{1,3}^{j-l}X_{l,1;j}^{1;4})-X_{j+1,1;j}^{1;4}}\\
&&\hspace{3cm}\times\left(\sum_{l=0}^{j}(-1)^{j-l}x_{1,3}^{j-l}x_{l,4}\right)\hspace{3cm}\mod \left<s_j\right>_{R_2\acute{}\,\acute{}}\\
&=&\frac{F_{j+1}(X_{1,1;j}^{1;4},\cdots,X_{j,1;j}^{1;4},x_{1,3}(\sum_{l=0}^{j}(-1)^{j-l}x_{1,3}^{j-l}X_{l,1;j}^{1;4}))-F_{j+1}(X_{1,1;j}^{1;4},\cdots,X_{j,1;j}^{1;4},X_{j+1,1;j}^{1;4})}{x_{1,3}-x_{1,1}}\\
&=&\frac{F_1(x_{1,3})-F_1(x_{1,1})}{x_{1,3}-x_{1,1}}+\frac{F_{j}(u_1,\cdots,u_j)-F_{j}(x_{1,4},\cdots,x_{j,4})}{x_{1,3}-x_{1,1}},
\end{eqnarray*}
where $\displaystyle u_k=\sum_{l=0}^{k}(-1)^{k-l}x_{1,3}^{k-l}X_{l,1;j}^{1;4}$.\\
The second term $\displaystyle\frac{F_{j}(u_1,u_2,\cdots,u_j)-F_{j}(x_{1,4},x_{2,4},\cdots,x_{j,4})}{x_{1,3}-x_{1,1}}$ equals
\begin{eqnarray*}
&&\frac{F_{j}(u_1,u_2,\cdots,u_j)-F_{j}(x_{1,4},u_2,\cdots,u_j)}{x_{1,3}-x_{1,1}}\\
&&+\frac{F_{j}(x_{1,4},u_2,u_3,\cdots,u_j)-F_{j}(x_{1,4},x_{2,4},u_3,\cdots,u_j)}{x_{1,3}-x_{1,1}}\\
&&+\cdots +\frac{F_{j}(x_{1,4},\cdots,x_{j-1,4},u_j)-F_{j}(x_{1,4},\cdots,x_{j-1,4},x_{j,4})}{x_{1,3}-x_{1,1}}.
\end{eqnarray*}
Since $u_k-x_{k,4}=(x_{1,1}-x_{1,3})\left(\sum_{l=0}^{k-1}(-1)^{k-1-l}x_{1,3}^{k-1-l}x_{l,4}\right)$ and $s_k=x_{k,2}-u_k$, the above polynomial equals
\begin{eqnarray*}
&&\frac{F_{j}(u_1,u_2,\cdots,u_j)-F_{j}(x_{1,4},u_2,\cdots,u_j)}{u_1-x_{1,4}}\\
&&+\frac{F_{j}(x_{1,4},u_2,u_3,\cdots,u_j)-F_{j}(x_{1,4},x_{2,4},u_3,\cdots,u_j)}{u_2-x_{2,4}}
\left(-\sum_{l=0}^{1}(-1)^{1-l}x_{1,3}^{1-l}x_{l,4}\right)\\
&&+\cdots +\frac{F_{j}(x_{1,4},\cdots,x_{j-1,4},u_j)-F_{j}(x_{1,4},\cdots,x_{j-1,4},x_{j,4})}{u_j-x_{j,4}}
\left(-\sum_{l=0}^{j-1}(-1)^{j-1-l}x_{1,3}^{j-1-l}x_{l,4}\right)\\
&\equiv&
\frac{F_{j}(x_{1,2},x_{2,2},\cdots,x_{j,2})-F_{j}(x_{1,4},x_{2,2},\cdots,x_{j,2})}{x_{1,2}-x_{1,4}}\\
&&+\frac{F_{j}(x_{1,4},x_{2,2},x_{3,2},\cdots,x_{j,2})-F_{j}(x_{1,4},x_{2,4},x_{3,2},\cdots,x_{j,2})}{x_{2,2}-x_{2,4}}
\left(-\sum_{l=0}^{1}(-1)^{1-l}x_{1,3}^{1-l}x_{l,4}\right)\\
&&+\cdots +\frac{F_{j}(x_{1,4},\cdots,x_{j-1,4},x_{j,2})-F_{j}(x_{1,4},\cdots,x_{j-1,4},x_{j,4})}{x_{j,2}-x_{j,4}}
\left(-\sum_{l=0}^{j-1}(-1)^{j-1-l}x_{1,3}^{j-1-l}x_{l,4}\right)\\
&&\hspace{7cm}\mod \left<s_1,s_2,\cdots,s_j\right>_{R_2\acute{}\,\acute{}}\\
&=&-L_{1,j}^{2;4}-L_{2,j}^{2;4}\left(\sum_{l=0}^{1}(-1)^{1-l}x_{1,3}^{1-l}x_{l,4}\right)-\cdots -L_{j,j}^{2;4}\left(\sum_{l=0}^{j-1}(-1)^{j-1-l}x_{1,3}^{j-1-l}x_{l,4}\right).
\end{eqnarray*}
Hence using Theorem \ref{reg-eq} and Proposition \ref{equiv}, the matrix factorization (\ref{total-mat3}) is isomorphic to
\begin{eqnarray*}
&&K\left(\left(
\begin{array}{c}
\ast\\[.5em]
\vdots\\[.5em]
\ast\\[.6em]
\frac{b}{\beta}
\end{array}
\right)
;
\left(
\begin{array}{c}
s_1\\[.5em]
\vdots\\[.5em]
s_j\\[.5em]
L_{1,1}^{3;1}-\sum_{k=1}^{j}L_{k,j}^{2;4}\left(\sum_{l=0}^{k-1}(-1)^{k-1-l}x_{1,3}^{k-1-l}x_{l,4}\right)
\end{array}
\right)
\right)_{R_2\acute{},\acute{}}\{1-n\}\left<1\right>\\
&\simeq&
K\left(\left(
\begin{array}{c}
L_{1,j}^{2;4}\\[.5em]
L_{2,j}^{2;4}\\[.8em]
\vdots\\[.5em]
L_{j,j}^{2;4}\\[.8em]
x_{1,3}-x_{1,1}
\end{array}
\right)
;
\left(
\begin{array}{c}
x_{1,2}-x_{1,4}+x_{1,3}-x_{1,1}\\[.5em]
x_{2,2}-x_{2,4}+(x_{1,3}-x_{1,1})\left(\sum_{l=0}^{1}(-1)^{1-l}x_{1,3}^{1-l}x_{l,4}\right)\\[.5em]
\vdots\\[.5em]
x_{j,2}-x_{j,4}+(x_{1,3}-x_{1,1})\left(\sum_{l=0}^{j-1}(-1)^{j-1-l}x_{1,3}^{j-1-l}x_{l,4}\right)\\[.5em]
L_{1,1}^{3;1}-\sum_{k=1}^{j}L_{k,j}^{2;4}\left(\sum_{l=0}^{k-1}(-1)^{k-1-l}x_{1,3}^{k-1-l}x_{l,4}\right)
\end{array}
\right)
\right)_{R_2\acute{},\acute{}}\{1-n\}\left<1\right>\\
&\simeq&
K\left(\left(
\begin{array}{c}
L_{1,j}^{2;4}\\[.5em]
\vdots\\[.5em]
L_{j,j}^{2;4}\\[.5em]
x_{1,3}-x_{1,1}
\end{array}
\right)
;
\left(
\begin{array}{c}
x_{1,2}-x_{1,4}\\[.8em]
\vdots\\[.5em]
x_{j,2}-x_{j,4}\\[.8em]
L_{1,1}^{3;1}
\end{array}
\right)
\right)_{R_2\acute{},\acute{}}\{1-n\}\left<1\right>\\
&\simeq&
K\left(\left(
\begin{array}{c}
L_{1,j}^{2;4}\\[.5em]
\vdots\\[.5em]
L_{j,j}^{2;4}\\[.5em]
L_{1,1}^{3;1}
\end{array}
\right)
;
\left(
\begin{array}{c}
x_{1,2}-x_{1,4}\\[.8em]
\vdots\\[.5em]
x_{j,2}-x_{j,4}\\[.8em]
x_{1,3}-x_{1,1}
\end{array}
\right)
\right)_{R_2\acute{},\acute{}}\simeq \c\acute{}\left( \input{figsquare1j-rev-mf}\right)_n.
\end{eqnarray*}
\end{proof}
\begin{cor}\label{cor-square}
\begin{eqnarray*}
&&\c\acute{}\left( \input{figsquare1j--k--k+1j-k--k--1j-mf2}\right)_n\\[-0.5em]
&\simeq&
\bigoplus_{i=0}^{(j_1-j_2)(j_2-1)}\left(\c\acute{}\left( \input{figsquare1j1--j1+1--1j1-mf}\right)_n\{-(j_1-j_2)(j_2-1)+2i\}\right)^{p(i,j_1-1;j_2-1)}\\[-0.5em]
&&\hspace{2cm}\bigoplus_{i=0}^{(j_1-j_2-1)j_2}\left(\c\acute{}\left( \input{figsquare1j1-mf}\right)_n\{-(j_1-j_2-1)j_2+2i\}\right)^{p(i,j_1-1;j_2)}
\end{eqnarray*}
\end{cor}
\begin{proof}
We consider the following matrix factorization
\begin{equation}\label{gene-cor}
\c\acute{}\left(\input{fig-gene-cor}\right)_n.
\end{equation}
Using Proposition \ref{mat-equiv3} (2) and Proposition \ref{mat-equiv2} (1), the matrix factorization is isomorphic to
\begin{eqnarray}
\nonumber&&\c\acute{}\left( \input{figsquare1j--k--k+1j-k--k--1j-mf2}\right)_n\bigoplus_{i=1}^{j_2-1}\c\acute{}\left(\input{figsquare1j-bubble-mf}\right)_n\{2i-j_2\}\\[-0.5em]
\nonumber
&\simeq&\c\acute{}\left( \input{figsquare1j--k--k+1j-k--k--1j-mf2}\right)_n\\[-0.5em]
\label{part-mat1}&&{}\bigoplus_{i=1}^{j_2-1}\left(\bigoplus_{k=0}^{j_1j_2-j_2^2}\left(\c\acute{}\left(\input{figsquare1j1-mf}\right)_n\{-j_1j_2+j_2^2+2k\}\right)^{p(k,j_1;j_2)}\right)\{2i-j_2\}.
\end{eqnarray}
On the other hand, using Proposition \ref{mat-equiv1}, Proposition \ref{mat-equiv2} (1) and Proposition \ref{mat-equiv3} (2), the matrix factorization (\ref{gene-cor}) is isomorphic to
\begin{eqnarray}
\nonumber&&\c\acute{}\left(\input{figsquare1j1--1--2j1-1-bubble-mf}\right)_n\\[-0.5em]
\nonumber&\simeq&\bigoplus_{i=0}^{(j_1-j_2)(j_2-1)}\left(\c\acute{}\left(\input{figsquare1j1--1--2j1-1--1--1j1-mf}\right)_n\{-(j_1-j_2)(j_2-1)+2i\}\right)^{p(i,j_1-1;j_2-1)}\\[-0.5em]
\nonumber&\simeq&
\bigoplus_{i=0}^{(j_1-j_2)(j_2-1)}\left(\c\acute{}\left( \input{figsquare1j1--j1+1--1j1-mf}\right)_n\{-(j_1-j_2)(j_2-1)+2i\}\right)^{p(i,j_1-1;j_2-1)}\\[-0.5em]
\label{part-mat2}&&\hspace{1cm}\bigoplus_{i=0}^{(j_1-j_2)(j_2-1)}\left(\left(\bigoplus_{k=1}^{j_1-1}\c\acute{}\left( \input{figsquare1j1-mf}\right)_n\{2k-j_1\}\right)\{-(j_1-j_2)(j_2-1)+2i\}\right)^{p(i,j_1-1;j_2-1)}.
\end{eqnarray}
The $\Z$-grading shift of the matrix factorization (\ref{part-mat1}) is derived from $\displaystyle [j_2-1]\left[j_1\atop j_2\right] $ and one of the matrix factorization (\ref{part-mat2}) is also derived from $\displaystyle [j_1-1]\left[j_1-1\atop j_2-1\right] $.
Since $\displaystyle [j_1-1]\left[j_1-1\atop j_2-1\right]-[j_2-1]\left[j_1\atop j_2\right]=\left[j_1-1\atop j_2\right] $, we obtain Corollary \ref{cor-square}.
\end{proof}

\begin{conj}[For Strategy (S1)]
There are equivalences in the category $\HMF^{gr}$ corresponding to the remaining MOY relations.
\end{conj}

For crossings with coloring $(1,2)$ and $(2,1)$, the author categorified MOY polynomial using the category $\kom(\HMF)$ \cite{Yone2}.

%%%%%%%%%%%%%%%%%%%%%%%%%%%%%%%%%%%%%%%%%%%%%%%%%%%%%%%%%%%%%%%%%%%%%%%%%%%%%%%%%%%%%%%%%%%%%%%%%%%%%%%%
%
%
% Reference
%
%
%%%%%%%%%%%%%%%%%%%%%%%%%%%%%%%%%%%%%%%%%%%%%%%%%%%%%%%%%%%%%%%%%%%%%%%%%%%%%%%%%%%%%%%%%%%%%%%%%%%%%%%%


\begin{thebibliography}{references}
\bibitem{Bar}
D. Bar-Natan,
Khovanov's homology for tangles and cobordisms,
Geom. Topol.\textbf{9} (2005), 1443--1499 . 
\bibitem{Gepner}
D. Gepner, 
Fusion rings and geometry, 
Comm. Math. Phys. 141 (1991), no. 2, 381--411. 
\bibitem{Gornik}
B. Gornik,
Note on Khovanov link cohomology,
arXiv:math.QA/0402266.
\bibitem{K1}
M. Khovanov,
A categorification of the Jones polynomial,
Duke Math. J. \textbf{101} (2000), no. 3, 359-426.
\bibitem{KR1}
M. Khovanov, L. Rozansky,
Matrix factorizations and link homology,
arXiv:math.QA/0401268.
\bibitem{KR2}
M. Khovanov, L. Rozansky,
Matrix factorizations and link homology II,
arXiv:math.QA/0505056.
\bibitem{KR3}
M. Khovanov, L. Rozansky,
Virtual crossings, convolutions and a categorification of the SO(2N) Kauffman polynomial,
arXiv:math.QA/0701333
\bibitem{MOY}
H. Murakami, T. Ohtsuki, S. Yamada,
Homfly polynomial via an invariant of colored plane graphs,
Enseign. Math. (2) \textbf{44} (1998), no. 3-4, 325--360.    
\bibitem{Ras}
J. Rasmussen,
Some differentials on Khovanov-Rozansky homology,
arXiv:math.GT/0607544.
\bibitem{Yone1}
Y. Yonezawa,
Matrix factorizations and double line in $\mathfrak{sl}_n$ quantum link invariant, arXiv:math.GT/0703779.
\bibitem{Yone2}
Y. Yonezawa,
A categorification of MOY polynomial with coloring 1 and 2.
\bibitem{Wu}
H. Wu,
On the quantum filtration of the Khovanov-Rozansky cohomology,
arXiv:math.GT/0612406.
\bibitem{Wu2}
H. Wu,
Matrix factorizations and colored MOY graphs,
arXiv:0803.2071.
\end{thebibliography}
\end{document}